\documentclass[reqno]{amsart}
\usepackage[english]{babel}
\usepackage{amsthm}
\usepackage{amsmath}

\usepackage{amssymb,amsthm,amsfonts}
\usepackage[T1]{fontenc}
\usepackage[latin1]{inputenc}
\usepackage[dvips]{graphicx}
\usepackage{fancyhdr}
\usepackage{color}
\usepackage{textcomp}
\usepackage{verbatim}

\newcommand{\bm}[1]{\text{\boldmath $#1$}}

\allowdisplaybreaks[3]
\renewcommand{\div}{\mathop{\rm div}\nolimits}
\newcommand{\inte}{\mathop{\rm int}\nolimits}
\newcommand{\A}{\mathcal{A}}
\newcommand{\B}{\mathcal{B}}
\newcommand{\R}{\mathbb{R}}

\newcommand{\Chi}{\boldsymbol{\chi}}

\newenvironment{pr}{\begin{proof}[\textbf{Proof}]}{\end{proof}}
\newtheorem{teor}{Theorem}[section]
\newtheorem{lemm}[teor]{Lemma}

\newtheorem{propos}[teor]{Proposition}

\newtheorem{obs}[teor]{Remark}

\numberwithin{equation}{section}
\providecommand{\abs}[1]{\left\lvert#1\right\rvert}
\providecommand{\norm}[1]{\left\lVert#1\right\rVert}

\begin{document}
\title{A dissipative model for hydrogen storage: \\existence and regularity results}

\author{Elisabetta Chiodaroli}

\address{Elisabetta Chiodaroli,
Institut f\"ur Mathematik, Universit\"at Z\"urich,
Winterthurerstrasse 190, CH-8057 Z\"urich, Switzerland}

\email{elisabetta.chiodaroli@math.uzh.ch}


\begin{abstract} We prove global existence of a solution to an
initial and boundary value problem for a highly nonlinear PDE
system. The problem arises from a thermomechanical dissipative model
describing hydrogen storage by use of metal hydrides. In order to
treat the model from an analytical point of view, we formulate it as
a phase transition phenomenon thanks to the introduction of a
suitable phase variable. Continuum mechanics laws lead to an
evolutionary problem involving three state variables: the
temperature, the phase parameter and the pressure. The problem thus
consists of three coupled partial differential equations combined
with initial and boun\-dary conditions. Existence and regularity of
the solutions are here investigated by means of a time
discretization-a priori estimates-passage to the limit procedure
joined with compactness and monotonicity arguments.
\end{abstract}

\maketitle

\pagestyle{myheadings}

\section{Introduction}
The paper deals with a thermo-mechanical model describing hydrogen
sto\-ra\-ge in terms of metal hydrides. Hydrogen storage basically
implies a reduction in the enormous volume of hydrogen gas; $1$ kg
of hydrogen at ambient temperature and atmospheric pressure has a
volume of $11$ $m^3$. There are basically six methods in order to
store hydrogen reversibly with a high volumetric and gravimetric
density: hydrogen gas, liquid hydrogen, physisorption, complex,
chemical hydrides, metal hydrides. This latter technique exploits
the possibility of many metals to absorb hydrogen: such metals and
alloys are able to react spontaneously with hydrogen and they can
store a large amount of it. These materials, either a defined
compound or a solid solution, are designed as metal hydrides: their
use in the hydrogen storage is of interest in terms of safety,
global yield and long-time storage. Energetic and industrial
applications warrant this interest. Indeed, hydrogen is foreseen to
be a clean and efficient energy carrier for the future (cf.
\cite{la}). Nowadays, energetic needs are mainly covered by fossil
energies leading to pollutant emissions mostly responsible for
global warming. Hydrogen stands among the best solutions to the
shortage of fossil energies and to the greenhouse effect, in
particular for energy transportation.

Our research moves in the direction of providing a predictive
theo\-ry for the storage of hydrogen by use of metal hydrides (cf.
\cite{la}). To this aim our analysis refers to the thermo-mechanical
model introduced by E. Bonetti, M. Fr\'emond and C. Lexcellent in
\cite{bfl} and it complements their results. Following the usual
approach of Thermodynamics, Bonetti, Fr\'emond and Lexcellent have
derived the governing equilibrium equations. The analytical
formulation of the problem they obtained was new and, to the best of
our knowledge, no other related results can be found in the
literature. Our work is attempted to go one step further both in
modeling and analytical aspects of hydrogen storage by use of metal
hydrides (cf. \cite{chio}). The problem we will deal with results
from a phase tran\-si\-tion model and is formulated via the
classical principles of Con\-ti\-nuum Me\-cha\-nics. The related
strongly nonlinear PDE system has been investigated, from the point
of view of exi\-stence and regularity of the solutions.

In order to get acquainted with the
phenomenon, we recall some of its basic features. Some metals are able to absorb hydrogen atoms and
combine with them to form solid solutions. We assume the existence
of two solid solutions: the $\alpha$-phase and the $\beta$-phase.
The presence of one phase with respect to the other depends on the
pressure of the hydrogen. To provide a good mechanical model to be
analytically treated from the point of view of existence and
regularity of solutions, it seems useful to exploit the theory of
phase transitions. We choose the volume fraction of one of the
phases as a state quantity and we denote it by $\chi$. Hence, $\chi$
satisfies the relation
\begin{equation} \label{internalconstraint}
\chi\in [0,1],
\end{equation}
and, assuming that no voids appear in the mixture, the volume
fraction of the other phase is simply given by $1-\chi$. More
precisely, if $\chi=1$ we have the $\alpha$-phase, if $\chi=0$ we
have the $\beta$-phase, and if $\chi \in (0,1)$ both phases are
present in suitable proportions. The state variables of the model
are the absolute temperature $\theta$, the hydrogen pressure $p$,
and the phase parameter $\chi$ along with its gradient $\nabla\chi$
accounting for local interactions between the different phases.
Constitutive relations for the state quantities will be chosen in
such a way that the principles of Thermodynamics are satisfied.
Finally, to describe the thermo-mechanical evolution of the system,
a pseudo-potential of dissipation will be considered. The
constitutive equations will be substantially recovered as in
\cite{bfl}. Nevertheless, some improvements will be introduced in
the model. Indeed, in \cite{bfl} the PDE system is written by
neglecting dissipative effects and microscopic movements in the
power of interior forces. Our approach will be different, as we aim
to derive the complete dissipative model accounting for microscopic
velocities and diffusive phase transformations. Moreover, \cite{bfl}
is concerned with the study of a weaker formulation of the system:
in this framework the authors have been able to prove a global
existence result holding for $n=3$. We complement their results.
Indeed, in our framework we can prove a more general existence
result in the three-dimensional setting. Moreover, by refining the
assumptions on the data, we are able to show further regularities
for the solutions to our problem. Some related analytical results
can be found for models describing irreversible phase transition
phenomena. In particular, we refer to \cite{lss} and \cite{lu}:
these papers deal with nonlinear systems of PDEs go\-ver\-ning the
evolution of two unknown fields ($\theta$ and $\chi$) and prove some
existence results. We refer also to \cite{bo} and \cite{bon}, both
concerning the analysis of a dissipative Fr\'emond model for shape
memory alloys. The problem therein is quite similar to our both for
modeling aspects and analytical investigation: many techniques and
analytical tools from \cite{bo} and \cite{bon} have inspired our
proofs.

Now, let us introduce the complete thermo-mechanical model
describing hydrogen storage by use of metal hydrides and including
dissipative effects and miscroscopic velocities in the constitutive
equations as well as microscopic forces in the principle of virtual
power. Then, we state an initial and boundary value problem for the
obtained model.

We write the model in terms of the state variables $\theta$, $p$ and
$\chi$. In order to achieve a precise description of the phenomenon,
we take into account both the equilibrium and the evolution of the
system, which are characterized by two functionals: the free energy
$\Psi$, defined on state variables, and the pseudo-potential of
dissipation $\Phi$, defined on dissipative variables. Constitutive
relations for the involved thermo-mechanical quantities are chosen
in accordance with the principles of Thermodynamics.

We preface our discussion by defining the hydrogen density $\rho_H$
and the total density $\rho$ (in what follows we take $\rho=1$).
Thus, letting $c_H$ be the capacity of the hydrogen, we have that
$c_H=\rho_H(\rho-\rho_H)^{-1}$, i.e.,
$$\rho_H=\rho \frac{c_H}{1+c_H}=\frac{1}{\tau},$$
where $\tau$ stands for the hydrogen's specific volume.

Now, we introduce the free energy function $\Psi(\theta, \tau, \chi,
\nabla\chi)$. By thermodynamical and duality arguments it follows
that the free energy is concave with respect to $\theta$, while we
assume that it is convex with respect to $\tau, \chi$, and
$\nabla\chi$. Thus, as $\Psi$ is convex with respect to $\tau$, we
can introduce the dual function $\Psi^\ast$ of $\Psi$ as follows:
$$\Psi^\ast(\theta, \zeta, \chi,\nabla\chi)=\sup_s\{\zeta s-\Psi(\theta, s, \chi,
\nabla\chi)\}.$$ Next, we define the pressure $p$ connected with the
hydrogen's specific volume $\tau$ by the following relation
\begin{equation} \label{definizionedip}
-p:=\frac{\partial\Psi}{\partial\tau}(\theta,\tau,\chi,\nabla\chi).
\end{equation}
This corresponds to setting $$\tau=\frac{\partial
\Psi^\ast}{\partial (-p)}(\theta,-p,\chi,\nabla\chi)=-\frac{\partial
\Psi^\ast}{\partial p}(\theta,-p,\chi,\nabla\chi).$$ Thus, if we
assume sufficient regularity for the functionals we get
$$\Psi^\ast(\theta, -p, \chi,
\nabla\chi)=-p\tau-\Psi(\theta, \tau, \chi, \nabla\chi).$$ Finally,
we deal with the \textit{enthalpy functional} $G(\theta, p, \chi,
\nabla\chi)$ defined in terms of the Legendre-Fenchel transformation
of $\Psi$ with respect to the specific volume $\tau$. More
precisely, the enthalpy $G(p)$ is defined by
$$G(p):=-\Psi^\ast(-p)$$ so that it results
\begin{equation}   \label{G}
G(\theta, p,\chi,\nabla\chi)=\Psi(\theta, \tau, \chi,
\nabla\chi)+p\tau.
\end{equation}
In particular, we recover that $G$ is concave with respect to $p$
and $\theta$, while it is convex with respect to $\chi$ and
$\nabla\chi$.

Then, we make precise the constitutive relations holding for the
entropy $s$, the specific volume $\tau$, and the internal energy $e$
(see \eqref{G}). We have
\begin{align}
&s=-\frac{\partial\Psi}{\partial\theta}=-\frac{\partial G}{\partial \theta}, \label{entropy}\\
&\tau=\frac{\partial G}{\partial p}, \label{volume}\\
&e=\Psi+\theta s=G-p\tau+\theta s. \label{internalenergy}
\end{align}
The fundamental balance laws of Continuum Mechanics, written in a
smooth bounded domain $\Omega\subset \R^n$ during a finite time
interval $[0,T]$ are: the momentum balance, the energy balance and
the mass balance. In the following we will denote the time
derivative of any function $f$ by the symbol $f_t$. At first, by the
principle of virtual power written for microscopic movements (i.e.,
neglecting any virtual macroscopic velocity), we recover an
equilibrium equation for the interior forces, which formally
correspond to the balance of the momentum seen as a microscopic
equilibrium equation. Namely, in absence of external actions, we
readily get
\begin{equation} \label{momentumbalance}
B-\div\bm{H}=0\quad \text{in }\Omega,
\end{equation}
coupled with the boundary condition
\begin{equation} \label{boundaryconditionforH}
\bm{H}\cdot\bm{n}=0\quad\text{on }\Gamma,
\end{equation}
where $\Gamma=\partial\Omega$. From
\eqref{momentumbalance}-\eqref{boundaryconditionforH} we can deduce
the mechanical meaning of the vector $\bm{H}$, which indicates a
work flux vector, while $B$ is a scalar quantity collecting
microscopic forces. In case when macroscopic deformations are
described by $-\tau_t$, we address the following energy balance
equation
\begin{equation} \label{energybalance}
e_t+\div\bm{q}=r+B\chi_t+\bm{H}\cdot\nabla\chi_t-p\tau_t \quad
\text{in }\Omega,
\end{equation}
where by $e$ we denote the internal energy as defined in
\eqref{internalenergy}; $\bm{q}$ represents the heat flux vector,
for which we will state later a suitable boundary condition. The
right hand side of \eqref{energybalance} accounts for heat sources
induced by mechanical and external actions. More precisely, $r$ is
an exterior heat source, while heat sources induced by microscopic
forces are collected by the terms involving the quantities $B$ and
$\bm{H}$. The presence of microscopic mechanically induced heat
sources in \eqref{energybalance} is justified by a generalization of
the principle of virtual power in which interior microscopic forces
and motions are also considered, as they are responsible for the
phase transition (cf. \cite{fre}). Finally \eqref{energybalance} is
complemented with a Neumann boundary condition
\begin{equation} \label{neumannconditionforq}
-\bm{q}\cdot\bm{n}=0 \quad \text{on }\Gamma,
\end{equation}
where $\bm{n}$ stands for the normal unit vector on the boundary
$\Gamma=\partial\Omega$. This corresponds to prescribe a null heat
flux through the boundary $\Gamma$. Here, and in the remainder of
the work, we assume small perturbations. Hence, letting the mass of
the hydrogen that is not in the solid solutions keep constant, the
hydrogen mass balance reads as follows
\begin{equation} \label{massbalance}
(\rho_H)_t+\div \bm{v}=0\quad\text{in }\Omega,
\end{equation}
where $\bm{v}$ is the hydrogen mass flux. Then we combine
\eqref{massbalance} with the following boundary condition
\begin{equation} \label{boundaryconditionforv}
-\bm{v}\cdot\bm{n}+\gamma p=0\quad \text{on }\Gamma, \quad \gamma>0,
\end{equation}
by which we require that the hydrogen flux through the boundary is
proportional to the difference between the exterior and the interior
pressure (here the exterior pressure is chosen equal to $0$).

In order to describe the thermo-mechanical evolution of the system
and to include dissipation in the model by following the approach by
Moreau (cf. \cite{mo}), we introduce a pseudo-potential of
dissipation $\Phi$ depending on $\chi_t$, $\nabla\chi_t$ and
$\nabla\theta$. We recall the properties of $\Phi$:
\begin{equation} \label{pseudo-potential}
\Phi\geq0, \qquad\Phi(\bm{0})=0, \qquad \Phi \quad \text{is convex.}
\end{equation}
By \eqref{pseudo-potential}, it turns out that the subdifferential
$\partial\Phi$ is a maximal monotone operator with
$\bm{0}\in\partial\Phi (\bm{0})$. In particular, it follows that
\begin{equation} \label{maximalmonotoneoperator}
\partial \Phi (\chi_t, \nabla \chi_t, \nabla \theta)\cdot
(\chi_t, \nabla\chi_t, ,\nabla\theta)\geq 0.
\end{equation}
Now, we are in the position to exhibit the constitutive relations
for $B$ and $\bm{H}$, in terms of $G$ and $\Phi$. Unlike Bonetti,
Fr\'emond and Lexcellent (see \cite{bfl}) we include dissipative
effects both in $B$ and in $\bm{H}$ and we prescribe them to be
given by the sum of a dissipative and a non-dissipative
contribution. Let us introduce a useful notation: nd in the apex is
used for pointing out nondissipative contributions, while d stands
for dissipative ones. Namely, we specify $B$ as
\begin{equation} \label{b}
B=B^{\text{nd}}+B^{\text{d}}=\frac{\partial G}{\partial
\chi}+\frac{\partial \Phi}{\partial \chi_t},
\end{equation}
and $\bm{H}$ as
\begin{equation} \label{H}
\bm{H}=\bm{H}^{\text{nd}}+\bm{H}^{\text{d}}=\frac{\partial
G}{\partial(\nabla\chi)}+\frac{\partial\Phi}{\partial(\nabla\chi_t)}.
\end{equation}
The heat flux vector $\bm{q}$ is assumed to fulfil the standard
Fourier law
\begin{equation} \label{fourierlaw}
\bm{q}=-k_0\nabla\theta,
\end{equation}
where $k_0>0$. Let us anticipate that, by a suitable choice of
$\Phi$, the heat flux can be expressed by use of the
pseudo-potential of dissipation. Finally, we set the following
relation for the hydrogen mass flux
\begin{equation} \label{v}
\bm{v}=-\lambda\nabla p,
\end{equation}
for $\lambda>0$ (take, e.g., $\lambda=1$).

Let us come to the functionals that describe the equilibrium and the
evolution of the system. As first we set the enthalpy functional $G$
as
\begin{equation} \label{enthalpy}
G(\theta, p, \chi, \tau)=a\log p+ b\chi(\log p-\log
p_e)-c_p\theta\log\theta+\frac{\delta}{2}\abs{\nabla\chi}^2+I_{[0,1]}(\chi),
\end{equation}
where $c_p>0$, $\delta>0$ and $I_{[0,1]}(\chi):=0$ if $\chi\in
[0,1]$ and $I_{[0,1]}(\chi):=+\infty$ otherwise. By $p_e$ in
\eqref{enthalpy} we denote the \textit{equilibrium} or
\textit{Plateau pressure} which is strongly temperature dependent.
In accordance with physical experience, we let $a>0$ and $b>0$
(take, e.g., $a=b=1$). Hence, experiments show that for $\theta$
sufficiently large the Van't Hoff law holds (see \cite{la}), i.e.,
\begin{equation} \label{Vanthoff}
\log p_e=-c_1\frac{1}{\theta}+c_2,
\end{equation}
where $c_1$ and $c_2$ are positive constants. However, as we have
already pointed out, the enthalpy $G$ has to be concave with respect
to the temperature, on the whole temperature interval. Thus we
complement relation \eqref{Vanthoff} by setting
\begin{equation}\label{functionh}
\log p_e=h(\theta),
\end{equation}
where $h$ is a sufficiently smooth function, e.g.
\begin{align}\label{hex}
&h(\theta)=-c_1\theta^{-1}+c_2\quad\text{for }\theta\text{
sufficiently large, say }\theta\geq\theta_{\ast}, \notag \\
&h(\theta)=\tilde{h}(\theta)\quad\text{for
}\theta_{\ast\ast}\leq\theta<\theta_{\ast},\notag\\
&h(\theta)=\text{constant}\quad\text{for }\theta<\theta_{\ast\ast},
\end{align}
with $\tilde{h}(\theta)$ suitably defined to yield $h\in C^2(\R)$.
Moreover, $\tilde{h}(\theta)$ and other values have to be chosen in
such a way that
\begin{equation} \label{necessarycondition}
\frac{\partial^2G}{\partial\theta^2}=-b\chi
h''(\theta)-\frac{c_p}{\theta}<0.
\end{equation}
Indeed, \eqref{necessarycondition} ensures that $G$ is concave with
respect to temperature, which is necessary to get the thermodynamic
consistency of the model. As it will be clear in what follows,
\eqref{necessarycondition} is in direct relationship with some
assumptions on $h$ concerning the analytical solvability of the
resulting heat equation.

In order to derive the model in a thermodynamic frame we have to
fulfil the second law of Thermodynamics in the form of the
Clausius-Duhem ine\-qua\-li\-ty. To this aim, it is convenient to
introduce the heat flux vector $\bm{q}$ formally as a dissipative
quantity defined by the pseudopotential of dissipation $\Phi$. Thus,
we set
\begin{equation} \label{phi}
\Phi=\frac{\mu}{2}
\abs{\chi_t}^2+\frac{\nu}{2}\abs{\nabla\chi_t}^2+\frac{k_0}{2\theta}\abs{\nabla\theta}^2,
\end{equation}
for $k_0>0$. Hence, we define a new dissipative quantity
\begin{equation} \label{Q}
\bm{Q}^d=-\frac{\partial\Phi}{\partial(\nabla\theta)}=-\frac{k_0}{\theta}\nabla\theta,
\end{equation}
so that letting
\begin{equation} \label{fourierlaw3}
\bm{q}=\theta\bm{Q}^d
\end{equation}
yields the Fourier relation \eqref{fourierlaw}. Now, we point out
that by use of the chain rule in \eqref{energybalance} and the above
constitutive relations \eqref{entropy}-\eqref{internalenergy},
\eqref{b}-\eqref{v}, \eqref{Q}-\eqref{fourierlaw3}, we can
equivalently rewrite the energy balance \eqref{energybalance} as
\begin{equation} \label{energybalancenew}
\theta \left( s_t+\div \bm{Q}^{\text{d}}-\frac{r}{\theta}
\right)=-\bm{Q}^{\text{d}}\cdot\nabla\theta+B^{\text{d}}\cdot\chi_t+\bm{H}^{\text{d}}\cdot\nabla\chi_t.
\end{equation}
Thus, after observing that by \eqref{b}, \eqref{H} and \eqref{Q} the
right hand side of \eqref{energybalancenew} corresponds to
\begin{equation} \label{righthandside}
\left(\frac{\partial\Phi}{\partial\chi_t},
\frac{\partial\Phi}{\partial\nabla\chi_t},\frac{\partial\Phi}{\partial\nabla\theta}\right)\cdot(\chi_t,\nabla\chi_t,\nabla\theta),
\end{equation}
we get, thanks to \eqref{maximalmonotoneoperator},
\begin{equation} \label{righthandside2}
-\bm{Q}^{\text{d}}\cdot\nabla\theta+B^{\text{d}}\cdot\chi_t+\bm{H}^{\text{d}}\cdot\nabla\chi_t\geq0.
\end{equation}
Finally, dividing \eqref{energybalancenew} by the absolute
temperature $\theta$ yields the Clausius-Duhem inequality ensuring
thermodynamic consistency, namely
\begin{equation} \label{thermodynamicconsistency}
s_t+\div\bm{Q}^{\text{d}}-\frac{r}{\theta}\geq 0.
\end{equation}
\begin{obs} \label{ossaddends}
Indeed, we have to remark that all the addends in
\eqref{righthandside2} turn out to be non-negative. In particular,
we have
\begin{align}
&-\bm{Q}^{\text{d}}\cdot\nabla\theta\geq0 \label{addend1} \\
&B^{\text{d}}\cdot\chi_t\geq0, \label{addend2} \\
&\bm{H}^{\text{d}}\cdot\nabla\chi_t\geq0. \label{addend3} \
\end{align}
Thanks to \eqref{addend1}-\eqref{addend2}, if we take $\nu=0$ in
\eqref{phi}, the model still fulfils the Clausius-Duhem inequality.
Indeed, letting $\nu=0$ in \eqref{phi} corresponds to neglecting
dissipative effects in $\bm{H}$ (see \eqref{H}): this choice would
lead to the model derived by Bonetti, Fr\'emond and Lexcellent in
\cite{bfl} which actually is thermodynamically consistent.
\end{obs}

Now, we can write the system of PDE's in terms of the unknowns, by
substituting in \eqref{momentumbalance}, \eqref{energybalance} and
\eqref{massbalance}, the constitutive equations of the model
specified by enthalpy and pseudo-potential of dissipation. In
particular, we write the vectors of microscopic forces \eqref{b} and
\eqref{H} on account of \eqref{enthalpy} and \eqref{phi}
\begin{align}
&B=b(\log p-h(\theta))+\partial I_{[0,1]}(\chi)+\mu \chi_t,
\label{vectorb} \\
&\bm{H}=\delta\nabla\chi+\nu\nabla\chi_t. \label{vectorH}
\end{align}
Let us note that $\partial I_{[0,1]}$ in \eqref{vectorb} stands for
the subdifferential of the indicator function of the convex
$[0,1]\subset\R$ and it is obtained as a ge\-ne\-ra\-li\-zed
de\-ri\-va\-ti\-ve with respect to $\chi$ of the non smooth function
$I_{[0,1]}$ in \eqref{enthalpy}. Actually, in our analysis we will
consider a more general maximal monotone graph $\beta$ in place of
$\partial I_{[0,1]}$, still with non negative and bounded domain.
Thus, the complete PDE system originating from
\eqref{energybalance}, \eqref{massbalance} and
\eqref{momentumbalance} is written in $Q:=\Omega\times(0,T)$ as
follows:
\begin{align}
&(bh''(\theta)\theta\chi+c_p)\theta_t-k_0\Delta\theta=r+\mu\chi_t^2+\nu\abs{\nabla\chi_t}^2-b\theta
h'(\theta)\chi_t, \label{PDE1} \\
&\mu\chi_t-\nu\Delta\chi_t-\delta\Delta\chi+\xi= -b(-h(\theta)+\log p), \label{PDE2} \\
&\xi\in \beta(\chi), \label{internalconstraintnew} \\
&\left(\frac{p}{a+b\chi}\right)_t-\lambda\Delta p=0. \label{PDE3}
\end{align}
Then \eqref{PDE1}-\eqref{PDE3} are combined with initial conditions
\begin{align}
&\theta(0)=\theta_0, \label{thetainitial} \\
&\chi(0)=\chi_0, \label{chiinitial} \\
&p(0)=p_0, \label{pinitial}
\end{align}
and the natural boundary conditions ($\partial_n$ is the normal
derivative operator on the boundary $\Gamma$)
\begin{align}
&k_0\partial_n\theta=0, \label{thetaboundary} \\
&\nu \left(\partial_n\chi\right)_t+\delta\partial_n\chi=0, \label{chiboundary} \\
&\lambda\partial_n p+\gamma p=0, \label{pboundary}
\end{align}
on $\Gamma\times(0,T)$. Observe that the unusual boundary condition
\eqref{chiboundary} results from the position \eqref{vectorH} and,
however, it is perfectly equivalent to the standard boundary
condition $\partial_n\chi=0$ whenever the compatibility condition
$\partial_n\chi_0=0$ holds true on $\Gamma$, as in our approach (cf.
\eqref{chi0}).

Let us point out to the reader that we deal with a slightly modified
version of the system \eqref{PDE1}-\eqref{PDE3}, as in equation
\eqref{PDE1} we consider $\nu=0$ and neglect exterior heat sources
($r=0$). From a mechanical perspective, letting $\nu=0$ in
\eqref{PDE1} corresponds to require that dissipative effects on the
gradients of the phases are negligible in the power of interior
forces with respect to the other mechanically induced heat sources,
which is reasonable in the framework of the small perturbations
assumption. Thus in the sequel, by abuse of notation, we will refer
to \eqref{PDE1} but considering $\nu=0$ and $r=0$. Moreover, as our
existence theorem for weak solutions to problem
\eqref{PDE1}-\eqref{pboundary} does not ensure the pointwise
validity of \eqref{PDE2}, we have to consider an extension of
inclusion \eqref{internalconstraintnew} in the framework of a pair
of spaces in duality.

The above formulation of the problem of hydrogen storage complements
the model advanced in \cite{bfl}. In particular, the quadratic term
$\chi_t^2$ in \eqref{PDE1} is new, as well as $\Delta\chi_t$ in
\eqref{PDE2}. Let us note that the solvability of the resulting PDE
system, written as a phase-field problem, turns out to be
interesting from the analytical point of view. Indeed, the system
\eqref{PDE1}-\eqref{pboundary} is highly nonlinear and solving it
(in some suitable sense) requires non trivial analytical tools. This
is mainly due to the coupling of higher-order nonlinear
contributions involving the unknowns, a maximal monotone graph and a
quadratic dissipative term for the phase parameter, a logarithmic
term involving the pressure. More precisely, in the parabolic
equation \eqref{PDE1} the specific heat is a nonlinear function: to
ensure coerciveness, we need to prescribe a suitable assumption on
the function $h$. Dealing with the equation governing the evolution
of the pressure, we have to combine the regularity of the function
$\chi$ and the pressure $p$, mainly to control the nonlinear
evolution term. Finally, it is worth observing that the pressure has
a major role in the evolution of the phase through the presence of
its logarithm as a source in the corresponding evolution inclusion
(see \eqref{PDE2}). The logarithm is easily controlled for high
values of the pressure, whenever we are able to control $p$, but it
degenerates as $p\searrow0$. Thus, our proof present some ad hoc
estimates for \eqref{PDE3} to handle this nonlinearity in the phase
equation.

Here is the outline of the paper. In the next section, we introduce
an equivalent abstract formulation of the $n$-problem
\eqref{PDE1}-\eqref{PDE3}, on account of the initial and boundary
conditions \eqref{thetainitial}-\eqref{pboundary}. In particular, we
will focus our investigation on the three-dimensional problem as it
is more meaningful from the physical point of view. At first, under
suitable assumptions on the function $h$ in \eqref{functionh}, we
can state a global existence result (Theorem \ref{teorexistence})
holding for $\nu=0$ in \eqref{PDE1}. Finally, under refined assumptions, we establish some
further regularities for the state variables and the positivity of
the temperature (Theorems \ref{teorregularitypressure},
\ref{teorfurtherregularities} and \ref{teorpositivitytemperature}).

The existence result is proved in Sections 3 -- 4 -- 5 by exploiting
a semi-implicit time discretization scheme combined with an a priori
estimate-passage to the limit procedure. In section 4, we prove the
positivity of the temperature and further regularities for the state
variables, by performing suitable a priori estimates on the
solutions of the problem. In the Appendix, we present some results
we will refer to in our investigation.
\section{The abstract problem: main results}
Our analysis refers to an abstract version of the problem
\eqref{PDE1}-\eqref{pboundary}. We render the physical constants to
1 (i.e., $a=b=c_p=\lambda=k_0=\delta=\mu=\gamma=1$). Next, we start
by listing the main mathematical hypotheses of our work. At first,
we want to specify $\beta$: we introduce a convex, lower
semicontinuous and proper function
\begin{equation} \label{beta^}
\hat{\beta}: \mathbb{R}\rightarrow (-\infty, +\infty]
\end{equation}
satisfying the following property
\begin{equation} \label{beta^property2}
\inte D(\hat{\beta})=(0, \lambda_\beta),
\end{equation}
for some constant $\lambda_\beta >0$. By $D(\hat{\beta})$ we denote
the \textit{effective domain} of $\hat{\beta}$. In addition, we
assume that there exists $\chi_\ast\in(0,\lambda_\beta)$ such that
\begin{equation} \label{beta^property1}
0=\hat{\beta}(\chi_\ast)\leq \hat{\beta}(s)\quad\forall s\in
D(\hat{\beta}).
\end{equation}
Now, we can introduce the subdifferential of $\hat{\beta}$, i.e.,
\begin{equation} \label{beta}
\beta=\partial\hat{\beta};
\end{equation}
then, $\beta$ turns out to be a \textit{maximal monotone graph} in
$\R\times\R$ (we refer to \cite{bar} and \cite{bre} for basic
definitions and properties of maximal monotone operators). Let us
note that \eqref{beta^property1} implies $0\in\beta(\chi_\ast)$. In
the above positions $D=D(\beta)$ denotes the domain of the graph
$\beta$, i.e., $D=\{r\in\R:\beta(r)\neq\emptyset\}$; we recall that
$D\subseteqq D(\hat{\beta})$.
\begin{obs} \label{ossbeta}
The graph $\beta$ is introduced in order to generalize the graph
$\partial I_{[0,1]}$ in the physical model we derived; in particular
let us stress that the choice $\beta=\partial I_{[0,1]}$ would force
$\chi$ to attain values only in $[0,1]$ (see
\eqref{internalconstraint}). Nonetheless, more general classes of
constraints could be admissible for the phase variable than that
provided by $\xi\in\partial I_{[0,1]}(\chi)$. In particular, from a
ma\-the\-ma\-ti\-cal point of view, in this analysis we are able and
aim to deal with an arbitrary maximal monotone graph $\beta$ with
domain included in some bounded interval $[0,\lambda_\beta]$.
\end{obs}
According to Remark \ref{ossbeta}, in \eqref{internalconstraintnew}
we have considered $\beta(\chi)$ in place of $\partial
I_{[0,1]}(\chi)$. Thus, the system \eqref{PDE1}-\eqref{PDE3} can be
rewritten in $Q$ as follows (recall that $\nu=0$ and $r=0$ in
\eqref{PDE1}),
\begin{align}
&(h''(\theta)\theta\chi+1)\theta_t-\Delta\theta=\chi_t^2-\theta
h'(\theta)\chi_t, \label{PDE1new} \\
&\chi_t-\nu\Delta\chi_t-\Delta\chi+\xi= h(\theta)-\log p, \label{PDE2new} \\
&\xi\in \beta(\chi),  \label{PDE2'new} \\
&\left(\frac{p}{1+\chi}\right)_t-\Delta p=0. \label{PDE3new}
\end{align}
For the sake of clarity, before proceeding, we introduce some useful
notation. We set
$$ H:=L^2(\Omega), \quad V:=H^1(\Omega),$$
and identify $H$ with its dual space $H'$, so that
\begin{equation} \label{triplet}
V\hookrightarrow H \hookrightarrow V '
\end{equation}
with dense and compact injections. Let $(\cdot,\cdot)$ and
$\norm{\cdot}$ be the inner product and the corresponding norm in
$H$, and denote by $\langle\cdot,\cdot\rangle$ the duality pairing
between $V'$ and $V$. Hence, we introduce the following abstract
operators:
\begin{equation} \label{A}
\A:V\rightarrow V', \qquad \langle \A v, u\rangle=\int_\Omega \nabla
v \cdot \nabla u, \qquad u, v \in V
\end{equation}
\begin{equation} \label{B}
\B:V\rightarrow V', \qquad \langle\B v, u\rangle=\int_\Omega \nabla
v \cdot \nabla u +\int_\Gamma v u, \qquad u, v \in V.
\end{equation}
Then, to simplify notation, we set
$$W:=\{f\in H^2(\Omega): \partial_n f=0\quad\text{on}\quad \Gamma\}.$$
Clearly, we have that $W\subset V$ with compact embedding.  Whereas
the results we present refer also to lower dimensional space domains
we have restricted our in\-ve\-sti\-ga\-tion to the
three-dimensional framework as it is more meaningful from the
physical point of view. Thus, let
\begin{equation} \label{omega}
\Omega\subset\R^3\quad\text{and}\quad Q:=\Omega\times(0,T),
\end{equation}
where $T$ is a fixed final time. We assume $\Omega$ to be a smooth
bounded domain. We associate the functional $J_H$ and $J$ on $H$ and
$V$ to the function $\hat{\beta}$, as follows
\begin{align}
&J_H(v)=\int_\Omega \hat{\beta}(v) \quad \text{if }v\in H \text{ and
}\hat{\beta}(v)\in L^1(\Omega), \label{JH1} \\
&J_H(v)=+\infty \quad \text{if }v\in H \text{ and
}\hat{\beta}(v)\not\in
L^1(\Omega), \label{JH2} \\
&J(v)=J_H(v) \quad \text{if }v\in V. \label{JV}
\end{align}
As it is well known, $J_H$ and $J$ are convex and lower
semicontinuous on $H$ and $V$ respectively. Note that they are also
proper since $V$ contains all the constant functions. We denote by
$\partial_H J_H:H\rightarrow 2^{H}$ and by
$\partial_{V,V'}J:V\rightarrow 2^{V'}$ the corresponding
subdifferentials. We remind that
\begin{align}
\xi\in&\partial_{H} J_H(\chi) \quad \text{if and only if} \quad
\xi\in H, \notag \\
&\chi\in D(J_H),\quad\text{and} \quad
J_H(\chi)\leq(\xi,\chi-v)+J_H(v)\quad\forall v\in H
\label{subdifferentialofJH}\\
\xi\in&\partial_{V,V'} J(\chi) \quad \text{if and only if} \quad
\xi\in V', \notag \\
&\chi\in D(J),\quad\text{and} \quad
J(\chi)\leq\langle\xi,\chi-v\rangle+J(v)\quad\forall v\in V
\label{subdifferentialofJV}
\end{align}
where $D(\cdot)$ denotes, as above, the \textit{effective domain}
for functionals and multi\-valued operators. Note that $\partial_H
J_H:H\rightarrow 2^{H}$ and $\partial_{V,V'}J:V\rightarrow 2^{V'}$
are maximal monotone operators. Moreover, observe that for $\chi,
\xi\in H$ we have (see, e.g., \cite[Ex. 2.1.3, p. 21]{bre})
\begin{equation} \label{equivalencecondition}
\xi\in \partial_H J_H(\chi) \quad \text{if and only if} \quad
\xi\in\beta(\chi) \quad\text{almost everywhere in }\Omega.
\end{equation}
On the other hand, one can easily check that the inclusion
\begin{equation}\label{inclusionJHinJV}
\partial_H J_H(\chi)\subseteq H\cap \partial_{V,V'}J(\chi)\quad
\forall\chi \in V
\end{equation}
holds, just as a consequence of \eqref{JV} (compare
\eqref{subdifferentialofJV} with \eqref{subdifferentialofJH}).
Therefore, for $\chi\in V$ and $\xi\in H$, the condition
\begin{equation}\label{xiinbeta}
\xi\in\beta(\chi) \quad\text{almost everywhere in }\Omega
\end{equation}
implies
\begin{equation} \label{xiinsubdifferentialJV}
\xi\in\partial_{V,V'}J(\chi).
\end{equation}
In addition, if we assume $\chi\in V$, $\xi\in H$ and
\eqref{xiinsubdifferentialJV}, then \eqref{xiinbeta} holds (cf.
\cite{bcg}). In particular, on account of \eqref{inclusionJHinJV},
we achieve the validity of the following equality
\begin{equation} \label{subdifferentialequality}
\partial_H J_H(\chi)= H\cap \partial_{V,V'}J(\chi)\quad
\forall\chi \in V.
\end{equation}
Nonetheless, we should observe that \eqref{subdifferentialequality}
is false for more general fun\-ctio\-nals, as one can verify by
referring to the example presented in \cite{bcg}.

Now, we are in the position to rewrite the PDE system
\eqref{PDE1new}-\eqref{PDE3new}, combined with initial and boundary
conditions \eqref{thetainitial}-\eqref{pboundary}, in the abstract
setting of the triplet $(V, H, V')$. We have to remark that we are
not able here to deal with a strong version of
\eqref{PDE1new}-\eqref{PDE3new}. In particular we cannot deal with
the natural extension of $\beta$, i.e. the subdifferential
$\partial_H J_H$. Hence, the variational inclusion governing the
dynamics of the phase is written in the abstract setting of the
$V'-V$ duality pairing. Nonetheless, even if it cannot be stated
a.e. in $Q$, it retains its physical consistence since it forces the
phase to attain only meaningful values. Indeed, if
$\xi(t)\in\partial_{V,V'}J(\chi(t))$ for almost any $t\in(0,T)$, we
have in particular that
\begin{equation} \label{chiinD}
\chi(t)\in D(\hat{\beta}) \quad\text{a.e. in }\Omega.
\end{equation}
For instance, if we take $\hat{\beta}=I_{[0,1]}$, the abstract
relation $\xi(t)\in\partial_{V,V'}J(\chi(t))$ for a.a. $t\in (0,T)$
implies that $\chi(t)\in[0,1]_V$\footnote{By $[0,1]_V$ we denote the
convex: $[0,1]_V:=\{v\in V: v\in[0,1]$ a.e. in $\Omega$\}.}, and
consequently $\chi\in [0,1]$ a.e. in $Q$. Thus the system is
rewritten, in $V'$ and a.e. in $(0,T)$, as follows:
\begin{align}
&e_t +
\A\theta=-h(\theta)\chi_t+\chi^2_t,\label{1}\\
&e=\theta- \chi (h(\theta)-\theta h'(\theta)),\label{1'}\\
&\chi_t+\A(\nu\chi_t+\chi) +\xi= h(\theta)-\log p,\label{2}\\
&\xi\in\partial_{V,V'}J(\chi),\label{2'}\\
&u_t+ \B p=0,\label{3}\\
&u=\frac{p}{1+\chi}.\label{defu}
\end{align}
Let us note that, by \eqref{1'} and \eqref{defu}, we have introduced
two new auxiliary variables: $e$ and $u$. The variable $e$ can be
expressed as a function $\psi$ of the variables $\theta$ and $\chi$:
\begin{equation} \label{e}
e=\psi(\theta,\chi):=\theta-\chi (h(\theta)-\theta h'(\theta)).
\end{equation}
On the other hand, the variable $u$ has a precise physical meaning:
it can be interpreted as the normalized hydrogen density. Thanks to
\eqref{thetainitial}-\eqref{pinitial} we are allowed to set the
following initial conditions for $e$ and $u$
\begin{align}
&e(0)=e_0:=\theta_0-\chi_0(h(\theta_0)-\theta_0
h'(\theta_0)),\label{einitial}\\
&u(0)=u_0:=\frac{p_0}{1+\chi_0}.\label{uinitial}
\end{align}
\begin{obs} \label{ossinclusion}
We point out that, whenever $\xi\in\partial_{V,V'}J(\chi)$ and
$\xi\in H$ a.e. in $(0,T)$ we have that actually $\xi\in\partial_H
J_H (\chi)$ a.e. in (0,T), from which one can deduce that
$\xi\in\beta(\chi)$ a.e. in $Q$.
\end{obs}
Now, we set the hypothesis on the data prescribed in the first part
of our analysis. Concerning the Cauchy data in
\eqref{thetainitial}-\eqref{pinitial}, we assume that
\begin{align}
&\theta_0 \in V, \label{teta0}\\
&\chi_0 \in W\cap D(J_H)\label{chi0}\\
&p_0 \in V, \quad\log p_0 \in L^1(\Omega).\label{p0}
\end{align}
Note that \eqref{p0} yields $p_0>0$ a.e. in $\Omega$. Moreover,
\eqref{chi0} implies that $\chi_0\in D(\hat{\beta})$ a.e. in
$\Omega$ and in particular:
\begin{equation} \label{boundchi0}
0\leq\chi_0\leq\lambda_\beta\quad \text{a.e. in }\Omega,
\end{equation}
where the value $\lambda_\beta$ is introduced in
\eqref{beta^property2}. Finally, from \eqref{chi0}-\eqref{boundchi0}
we can deduce that
\begin{equation}\label{u0}
u_0\in V.
\end{equation}
In fact, the following estimate holds
\begin{equation} \label{estimateforu0}
\norm{u_0}^2_V\leq C\left(\norm{p_0}^2_V+\norm{p_0}^2_V\norm{\nabla
\chi_0}^2_V\right)\leq C,
\end{equation}
where we have exploited the continuous embedding $V\subset
L^4(\Omega)$.

Hence, we ask for a suitable regularity of the thermal expansion
coefficient $h(\theta)$, in agreement with the assumptions leading
to the physical consistence of the model (see
\eqref{necessarycondition}). We require
\begin{align}
&h \in W^{2,\infty}(\R)\cap C^2(\R),\label{h1}\\
&\norm{h}_{W^{2,\infty}(\R)}+\abs{h'(\zeta)\zeta}\leq c_h, \quad
\abs{h''(\zeta)\zeta}\leq c'_h, \quad \forall \zeta\in \R,\label{h2}
\end{align}
for some positive constants $c_h, c'_h$. In addition, let $c_s>0$
such that (recall \eqref{chiinD} and \eqref{h2})
\begin{equation} \label{h3}
1+\eta h''(\zeta)\zeta\geq c_s>0, \quad \forall \eta\in
D(\hat{\beta}), \quad\forall\zeta\in\R.
\end{equation}
This correspond to assume that the product $c'_h\lambda_\beta$ is
small with respect to $1$. The hypotheses we made on $h$ allow us to
infer that $\psi$ (see \eqref{e}) is a bi-lipschitz function of the
variable $\theta$. Indeed the following bounds hold (see
\eqref{h2}-\eqref{h3})
\begin{align}
&0<c_s\leq\partial_1 \psi\leq 1+c_e, \label{derivatapsi1}  \\
&\abs{\partial_2 \psi}\leq c_h, \label{derivatapsi2}
\end{align}
for a positive constant $c_e$ depending on $c_h'$, where
$\partial_1\psi$, $\partial_2\psi$ denote the partial derivatives of
$\psi$ with respect to first and second variable, respectively. Let
us recall that for the initial datum $e_0$ we have:
$e_0=\psi(\theta_0, \chi_0)$ (see \eqref{einitial}). This and the
regularity of $\theta_0$ and $\chi_0$ in \eqref{teta0}-\eqref{chi0},
along with the above properties of $\psi$, easily yield
\begin{equation}\label{e0}
e_0\in V.
\end{equation}
By use of a semi-implicit time discretization scheme combined with
the a priori estimate-passage to the limit procedure, we can prove
the fol\-lowing re\-la\-ted global existence result.
\begin{teor}[\textbf{Existence}] \label{teorexistence}
Let \eqref{teta0}-\eqref{p0}, \eqref{u0}, \eqref{e0} and
\eqref{h1}-\eqref{h3} hold. Then there exists a quintuple of
functions $(\theta,e, \chi, p, u)$ fulfilling
\begin{align}
&\theta\in L^\infty(0,T,H)\cap L^2(0,T,V), \label{tetaregularity} \\
&e\in W^{1,1}(0,T,V')\cap L^\infty(0,T,H)\cap L^2(0,T,V),\label{eregularity}\\
&\chi \in H^1(0,T,V)\cap L^\infty(0,T,W)\cap L^\infty(Q), \label{chiregularity} \\
&\hat{\beta}(\chi)\in L^\infty(0,T,L^1(\Omega)), \label{beta^regularity} \\
&p \in H^1(0,T,H)\cap L^\infty(0,T,V)\cap L^2(0,T, H^2(\Omega)), \label{pregularity} \\
&\log p \in L^\infty(0,T,L^1(\Omega))\cap L^2(0,T,V),\label{logpregularity}\\
&u \in H^1(0,T,H)\cap L^\infty(0,T,V), \label{uregularity}
\end{align}
and solving \eqref{1}-\eqref{defu} in $V'$, a.e. in $(0,T)$ along
with \eqref{chiinitial} and \eqref{einitial}-\eqref{uinitial}.
\end{teor}
The regularities obtained in Theorem \ref{teorexistence} allow the
equation \eqref{3} to make sense a.e. in $\Omega\times(0,T)$: the
same does not hold for \eqref{1} and \eqref{2} since the regularities
\textit{in space} are not strong enough. But, still dea\-ling with
the complete problem \eqref{1}-\eqref{defu}, we can establish some
further regularity results for the state variables and, in addition,
the positivity of the temperature $\theta$. As first, we can prove
some properties for the inverse of the pressure $p$. To this aim, we
need to make an additional assumption on the inverse of the initial
datum $p_0$. Namely, we need $1/p_0=p_0^{-1}\in H$. Then, the
following result holds.
\begin{teor}[\textbf{Regularity of the pressure}]\label{teorregularitypressure}
Under the same as\-sum\-ptions as in Theorem \ref{teorexistence}, if
\begin{equation} \label{p^-10}
\frac{1}{p_0}=p_0^{-1}\in H,
\end{equation}
then the following property holds
\begin{equation} \label{furtherpregularity}
 p^{-1}\in L^\infty (0,T,H)\cap L^2(0,T,V).
\end{equation}
\end{teor}
The next result is concerned with improved regularities for the time
derivatives of the state variables $\theta$ and $\chi$.
\begin{teor}[\textbf{Further Regularities}] \label{teorfurtherregularities}
Let \eqref{teta0}-\eqref{p0}, \eqref{u0}, \eqref{e0}, \eqref{p^-10}
and \eqref{h1}-\eqref{h3} hold. Let the quintuple $(\theta, e, \chi,
p, u)$ fulfill \eqref{tetaregularity}-\eqref{uregularity},
\eqref{furtherpregularity} and solve \eqref{1}-\eqref{defu} along
with \eqref{chiinitial} and \eqref{einitial}-\eqref{uinitial}. If
\begin{equation}\label{chi0new}
\chi_0\in D(\partial_{V,V'}J),
\end{equation}
then it holds that
\begin{align}
&\theta\in H^1(0,T,H)\cap
L^\infty(0,T,V),\label{furthertetaregularity}\\
&\chi\in W^{1,\infty}(0,T,V).\label{furtherchiregularity}
\end{align}
\end{teor}
Note that the improved regularity \eqref{furthertetaregularity} for
the temperature $\theta$ allow us to infer from \eqref{1'} that also
the time derivative of the internal variable $e$ belongs to
$L^2(0,T, H)$. Therefore, the solutions obtained in Theorem
\ref{teorfurtherregularities} fulfil equation \eqref{1} a.e. in
$\Omega\times(0,T)$. Finally, we can establish the positivity of the
temperature $\theta$, holding under suitable assumptions on the
logarithm of the initial datum $\theta_0$. Indeed, the following
theorem holds.
\begin{teor}[\textbf{Positivity of the temperature}] \label{teorpositivitytemperature}
Under the same assumptions as in Theorem
\ref{teorfurtherregularities}, if
\begin{equation} \label{logteta0}
\log \theta_0 \in L^1(\Omega),
\end{equation}
then the following properties
\begin{align}
&\log \theta \in L^\infty(0,T, L^1(\Omega))\cap L^2(0,T,V)
\label{logtetaregularity} \\
& \frac{\chi_t}{\sqrt{\theta}}\in L^2(0,T,H)
\label{furtherregularity}
\end{align}
hold.
\end{teor}

Clearly, \eqref{logtetaregularity} ensures the temperature $\theta$
to be positive a.e. in $Q$. Theorems \ref{teorregularitypressure},
\ref{teorfurtherregularities} and \ref{teorpositivitytemperature}
will be proved by performing some suitable a priori estimates on the
solutions of the problem whose existence is stated by Theorem
\ref{teorexistence}.

\section{Time discretization}
In order to achieve the proof of Theorem \ref{teorexistence} we
proceed as follows. First of all we establish a global existence and
uniqueness result for an approximating problem. Next, we perform
some a priori estimates that enable us to pass to the limit.

In this section, we approximate the system \eqref{1}-\eqref{defu} by
use of a semi-implicit time discretization scheme. Thus, letting $N$
be an arbitrary positive integer, we denote by $\tau:= T/N$ the time
step of our backward finite dif\-fe\-ren\-ces scheme. In the
forthcoming analysis, we will extensively use the following
notation. Let $(V^0,V^1,...,V^N)$ be a vector. Then, we denote by
$v_\tau$ and $\bar{v}_\tau$ two functions defined on the time
intervals $[0,T]$ and $(-\infty,T]$ which interpolate the values of
the vector piecewise linearly and backward constantly, respectively.
That is,
\begin{align}
v_\tau(0):&=V^0, \qquad \qquad v_\tau(t):=a_i(t)V^i+(1-a_i(t))V^{i-1} \label{tau} \\
\bar{v}_\tau(t):&=V^0 \text{ if } t\leq 0, \quad
\bar{v}_\tau(t):=V^i \text{ if } t\in((i-1)\tau, i\tau]
\label{bartau}\\
\intertext{where} a_i(t):&=\frac{(t-(i-1)\tau)}{\tau} \text{ if
}t\in((i-1)\tau, i\tau],
\end{align}
for $i=1,...,N$. Moreover, let us introduce the backward translation
operator $\mathcal{T_\tau}$ related to the time step $\tau$ by
setting
\begin{align}
&\mathcal{T_\tau} f(x,t):= f(x,t-\tau) \text{ for a.e. }
(x,t)\in\Omega\times (0,T),\notag \\
&\forall f:\Omega\times (-T,T)\rightarrow\R \text{
measurable}.\label{translation}
\end{align}
Next we regularize the initial datum for the internal variable $p$
by defining
\begin{equation} \label{p0tau}
p_{0\tau}=\begin{cases}
p_0&\text{ if }p_0\geq\tau\\
\tau &\text{ if }p_0<\tau.
\end{cases}
\end{equation}
Clearly, the regularized datum $p_{0\tau}$ is still in $V$. Let us
note that, thanks to \eqref{p0}, $\log p_{0\tau}\in L^1(\Omega)$ and
moreover
\begin{equation} \label{boundp0}
\norm{p_{0\tau}}_V+\norm{\log p_{0\tau}}_{L^1(\Omega)}\leq C,
\end{equation}
for some constant $C$ depending only on $\norm{p_0}_V$, $\norm{\log
p_0}_{L^1(\Omega)}$, $\abs{\Omega}$ and $T$. Besides, in view of
\eqref{p0}, \eqref{p0tau} ensures that
\begin{equation} \label{logp0tauinH}
\log p_{0\tau}\in H \quad\forall\tau
\end{equation}
as well. Let us note in advance that, in the approximating form, we
set the variational inclusion \eqref{2'} in $H$ by substituting the
abstract operator $\partial_{V,V'}J$ by the corresponding maximal
monotone graph $\partial_H J_H$ in $H$, provided we can prove some
regularity of the solutions. As a consequence we will be able to
solve the discrete variational inclusion a.e. in $\Omega$ (cf. also
Remark \ref{ossinclusion}). Then, the approximated problem can be
formulated as follows.\\
\textit{Problem $P_\tau$.} Find vectors
\begin{align}
&(\Theta^0, \Theta^1, ...., \Theta^N)\in V^{N+1},\label{thetadiscreti} \\
&(\Chi^0, \Chi^1, ...., \Chi^N)\in W^{N+1}, \label{chidiscreti} \\
&(P^0, P^1, ...., P^N)\in V^{N+1}, \label{pdiscreti} \\
\intertext{satisfying} \Theta^0&=\theta_0,\quad  \Chi^0=\chi_0,\quad
P^0=p_{0\tau},\label{condiniz}
\end{align}
and such that, by setting
\begin{align}
&E^0=\Theta^0-\Chi^0(h(\Theta^0)-\Theta^0 h'(\Theta^0)), \label{E0}\\
&U^0=\frac{P^0}{1+\Chi^0}, \label{U0}\\
\intertext{the following equations hold for $i=1,....,N$}
&\frac{E^i-E^{i-1}
}{\tau}+\A\Theta^i=-h(\Theta^{i-1})\frac{\Chi^i-\Chi^{i-1}}{\tau}+\left(\frac{\Chi^i-\Chi^{i-1}}{\tau}\right)^2
\quad\text{in}\quad V',\label{1approx}\\
&E^i=\Theta^i- \Chi^i(h(\Theta^i)-\Theta^i h'(\Theta^i)), \label{E}\\
&\frac{\Chi^i-\Chi^{i-1}}{\tau}+\nu\frac{\A\Chi^i-\A\Chi^{i-1}}{\tau}+
\A\Chi^i+\Xi^i=h(\Theta^{i-1})-\log P^{i-1}\quad \text{in}\quad
V',\label{2approx}\\
&\frac{U^i-U^{i-1}}{\tau}+\B P^i=0 \quad\text{in}\quad
V',\label{3approx}\\
&U^i=\frac{P^i}{1+\Chi^i},\label{U}\\
\intertext{for} &\Xi^i\in\partial_H J_H(\Chi^i) \label{2approx'}.
\end{align}
\begin{obs}
Since we have set the variational inclusion \eqref{2approx'} in $H$
and due to the regularity assumptions, we will see that equation
\eqref{2approx} turns out to make sense also in $H$.
\end{obs}
Let us note that \eqref{2approx'} implies that
$\Xi^i\in\beta(\Chi^i)$ a.e. in $\Omega$; in particular, since
$D(\hat{\beta})$ is included in some bounded interval
$[0,\lambda_\beta]$ (see \eqref{beta^property2} and
\eqref{boundchi0}), we can infer
\begin{equation} \label{chibounded}
0\leq\Chi^i\leq \lambda_\beta \quad\forall i=0,...,N,
\end{equation}
a.e. in $\Omega$. The conditions in \eqref{condiniz}, combined with
\eqref{E0} and \eqref{U0}, provide that (cf.
\eqref{einitial}-\eqref{uinitial})
\begin{align}
&E^0=e_0, \label{condinize}\\
&U^0=u_{0\tau}:=\frac{p_{0\tau}}{1+\chi_0}. \label{condinizu}
\end{align}
We can prove existence and uniqueness of a solution for the
approximating discrete problem $P_\tau$ at any step $\tau>0$.
Indeed, the following lemma holds.
\begin{lemm}[\textbf{Discrete well-posedness}] \label{lemma1}
Under the assumption \eqref{teta0}-\eqref{p0} and
\eqref{h1}-\eqref{h3}, for any $\tau>0$ the problem $P_\tau$ admits
a unique solution.
\end{lemm}
\begin{pr}
Owing to \eqref{condiniz}-\eqref{U0}, \eqref{teta0}-\eqref{p0} and
\eqref{p0tau}, we can restrict ourselves to prove that for any fixed
$\tau>0$ and for any $i\geq1$, the system
\eqref{1approx}-\eqref{2approx'} admits a unique solution. The main
idea is to proceed by induction on $i$. Indeed, we suppose to know
\begin{equation} \label{ipinduzione}
(\Theta^{i-1},\Chi^{i-1},P^{i-1}) \in V\times (W\cap D(J_H))\times
V.
\end{equation}
with
\begin{equation} \label{ipinduzionebis}
\log P^{i-1}\in H.
\end{equation}
(see \eqref{logp0tauinH}). We look for
\begin{equation} \label{tesiinduzione}
(\Theta^i,\Chi^i,P^i) \in V\times (W\cap D(J_H))\times V
\end{equation}
solving the resulting equations \eqref{1approx}-\eqref{2approx'} and
such that
\begin{equation} \label{logPiinH}
\log P^i\in H.
\end{equation}
We take first \eqref{2approx} and rewrite it as
\begin{equation}\label{2approxnew}
\frac{\Chi^i}{\tau}+\left(1+\frac{\nu}{\tau}\right)\A\Chi^i+\Xi^i=\frac{\Chi^{i-1}}{\tau}+\frac{\nu}{\tau}
\A\Chi^{i-1}+h(\Theta^{i-1})-\log P^{i-1}.
\end{equation}
Since we are assuming \eqref{h1}-\eqref{h2} and
\eqref{ipinduzione}-\eqref{ipinduzionebis} we can observe that the
right hand side, say $\mathcal{F}$, is known in $H$. Thus relation
\eqref{2approxnew} can be equivalently rewritten as
\begin{align}
&(\tau^{-1}Id+\mathcal{C}+\partial_H J_H)\Chi^i\ni\mathcal{F},\label{inclusionchi}\\
\intertext{where} &\mathcal{C}\Chi^i=(1+\nu\tau^{-1})\A\Chi^i.\notag
\end{align}
Actually, we would like to exploit the well known results on the
maximality of sums of monotone operators, holding under particular
regularity conditions. Namely, in this framework, we can invoke
Theorem \ref{teorsuriectivitymaximalmonotoneoperators} and get the
required existence of a function $\Chi^i\in W\cap D(J_H)$ fulfilling
\eqref{2approx}. On a second step, we take into consideration
\eqref{3approx}, where $\Chi^i$ is now the solution of
\eqref{2approx}. Equation \eqref{3approx} can be reformulated as
follows
\begin{equation}  \label{3approxnew}
\tau^{-1} \frac{P^i}{1+\Chi^i}+\B P^i=\tau^{-1}
\frac{P^{i-1}}{1+\Chi^{i-1}}
\end{equation}
The existence and uniqueness of a solution $P^i\in V$ follows
directly from the Lax-Milgram theorem, taking into account that
$\Chi^i$ satisfies \eqref{chibounded}. Before proceeding, we have to
check that \eqref{logPiinH} holds. To this aim, we refer to the
First a priori estimate, where we will prove that $\log U^i\in V$
(see \eqref{tesiinduzionebis}). In particular, this and the fact
that $\Chi^i$ obeys \eqref{chibounded} yield \eqref{logPiinH}.
Namely, owing to \eqref{U}, there is a constant $C$ such that
\begin{equation}
\norm{\log P^i}_H\leq C\left(\norm{\log
U^i}_H+\norm{\Chi^i}_H\right).
\end{equation}
Finally, letting in \eqref{1approx} ${\Chi^i}$ be the unique
solution of \eqref{2approx} and exploiting once more standard
results on maximal monotone operators we can find a unique function
$\Theta^i\in V$ solving the equation. Indeed we can rewrite
\eqref{1approx} more explicitly in terms of $\Theta^i$ as
\begin{align} \label{1approxnew}
&\tau^{-1} \Theta^i- \tau^{-1} \Chi^i (h(\Theta^i)-\Theta^i
h'(\Theta^i))+ \A \Theta^i \notag \\
&=  \tau^{-1}
\Theta^{i-1}(1+h'(\Theta^{i-1})\Chi^{i-1})-\tau^{-1}h(\Theta^{i-1})\Chi^i+
\left(\frac{\Chi^i-\Chi^{i-1}}{\tau}\right)^2.
\end{align}
As the right hand side, say $\mathcal{G}$, is known in $H$ the above
relation can be equivalently rewritten as
\begin{equation}  \label{inclusionteta}
(\tau^{-1} (Id+\mathcal{R}_i)+\A) \Theta^i\ni \mathcal{G},
\end{equation}
where $Id$ stands for the identity operator in $H$ and
$$\mathcal{R}_i(\Theta^i)=
\Chi^i(\Theta^ih'(\Theta^i)-h(\Theta^i))$$ maps $H$ into $H$. For
the sake of clarity, in our notation the subscript $i$ in
$\mathcal{R}_i$ is used for pointing out the dependence of
$\mathcal{R}_i$ on $\Chi^i$. Then, in view of assumptions
\eqref{h2}-\eqref{h3} it is not difficult to check that
$Id+\mathcal{R}_i$ is a Lipschitz continuous and strongly monotone
operator from $H$ into $H$, and consequently coercive. Since $\A$ is
a maximal monotone operator with domain $W$, the hypotheses of,
e.g., \cite[Cor 1.3, p. 48]{bar} hold and we get the required
existence and uniqueness of a solution $\Theta^i\in V$ to equation
\eqref{1approx}. This concludes our proof of Lemma \ref{lemma1},
since for any $i$, and any fixed $\tau>0$, the corresponding triple
$(\Theta^i,\Chi^i,P^i)$ solves the system
\eqref{1approx}-\eqref{2approx'}.
\end{pr}

\section{A priori estimates}
In this section, we aim to establish some a priori estimates on the
time-discrete solutions whose existence has been proved in Lemma
\ref{lemma1}. By virtue of Lemma \ref{lemma1} and owing to the
position \eqref{tau}-\eqref{translation}, we may introduce the
piecewise constant and linear in time functions $\bar{\theta}_\tau$,
$\bar{e}_\tau$, $\bar{\chi}_\tau$, $\bar{\xi}_\tau$, $\bar{p}_\tau$,
$\bar{u}_\tau$, $\theta_\tau$, $e_\tau$, $\chi_\tau$, $p_\tau$,
$u_\tau$ interpolating the corresponding values. Thanks to these
notations, the scheme \eqref{1approx}-\eqref{U} is restated as
follows in $V'$ and a.e. in $(0,T)$
\begin{align}
&\partial_t e_\tau+\A
\bar{\theta}_\tau=-h(\mathcal{T_\tau}\bar{\theta}_\tau)\partial_t\chi_\tau+(\partial_t\chi_\tau)^2,
\label{1tau} \\
&\bar{e}_\tau=\bar{\theta}_\tau-\bar{\chi}_\tau(h(\bar{\theta}_\tau)-\bar{\theta}_\tau
h'(\bar{\theta}_\tau)),\label{1tau'}\\
&\partial_t\chi_\tau+\nu\A(\partial_t\chi_\tau)+\A\bar{\chi}_\tau+\bar{\xi}_\tau
=h(\mathcal{T_\tau}\bar{\theta}_\tau)-\mathcal{T_\tau}(\log\bar{p}_\tau),
\label{2tau} \\
&\bar{\xi}_\tau\in \partial_{V,V'}J(\bar{\chi}_\tau), \label{2'tau}
\\
&\partial_t u_\tau+\B\bar{p}_\tau=0,\label{3tau}\\
&\bar{u}_\tau=\frac{\bar{p}_\tau}{1+\bar{\chi}_\tau}.
\label{defutau}
\end{align}
As concerns the inclusion \eqref{2approx} we could write it in $H$
in terms of the above introduced piecewise linear and constant
functions. Nonetheless, in order to perform a passage to the limit
procedure as $\tau\searrow0$, we have to set this inclusion in the
abstract framework of the $V-V'$ duality. Thus, instead of
$\partial_H J_H$ we have written the corresponding abstract operator
$\partial_{V,V'}J$ in \eqref{2'tau}. Obviously, by the regularity of
the solution vectors, the existence and uniqueness result we have
proved in the previous section can be extended to the abstract
framework of $V'$. In addition, let us observe that by construction
(cf. \eqref{tau}), $\chi_\tau$, $\theta_\tau$, $p_\tau$, $e_\tau$
and $u_\tau$ satisfy the natural Cauchy conditions (cf.
\eqref{p0tau}, \eqref{condiniz} and
\eqref{condinize}-\eqref{condinizu})
\begin{align}
&\chi_\tau(0)=\chi_0, \label{chitau0}\\
&\theta_\tau(0)=\theta_0,\label{tetatau0}\\
&p_\tau(0)=p_{0\tau},\label{ptau0} \\
&e_\tau(0)=e_0, \label{etau0}\\
&u_\tau(0)=u_{0\tau}. \label{utau0}
\end{align}
Hence, we are going to prove some estimates on the approximating
functions solving \eqref{1tau}-\eqref{defutau}; such estimates hold
at least for $\tau$ sufficiently small, but the involved constants
do not depend on $\tau$. Indeed, our aim is passing to the limit in
the above system as $\tau\searrow0$, by compactness or direct proof,
to get \eqref{1}-\eqref{defu} solved in a suitable sense. Let us
recall the trivial equality
\begin{equation} \label{eq}
2a(a-b)=a^2+(a-b)^2-b^2,\quad\forall a,b\in\R,
\end{equation}
which will be applied in the following estimates on the discrete
solutions.
\\
\\
\textbf{First a priori estimate.} We first test equation
\eqref{2approx} by $(\Chi^i-\chi_\ast)$, where the value $\chi_\ast$
is introduced in \eqref{beta^property1} and satisfies
$0\in\beta(\chi_\ast)$. Let us observe that, by monotonicity of the
operator $\partial_H J_H$, \eqref{2approx'} yields
\begin{equation} \label{1.1}
\int_\Omega\Xi^i(\Chi^i-\chi_\ast) \geq 0.
\end{equation}
Hence, exploiting the relation \eqref{eq}, by using \eqref{1.1}, we
get
\begin{align} \label{1.2}
&\frac{1}{2\tau}\norm{\Chi^i-\chi_\ast}^2_H+\frac{1}{2\tau}\norm{\Chi^i-\Chi^{i-1}}_H^2-\frac{1}{2\tau}\norm{\Chi^{i-1}-\chi_\ast}^2_H+
\left(1+\frac{\nu}{2\tau}\right)\norm{\nabla\Chi^i}^2_H
\notag\\
&+\frac{\nu}{2\tau}\norm{\nabla\Chi^i-\nabla\Chi^{i-1}}_H^2
-\frac{\nu}{2\tau}\norm{\nabla\Chi^{i-1}}^2_H\notag\\
&\leq\int_\Omega h (\Theta^{i-1})(\Chi^i-\chi_\ast)-\int_\Omega \log
P^{i-1}(\Chi^i-\chi_\ast).
\end{align}
Then, we would like to formally test equation \eqref{3approx} by
$1-1/U^i$. Nevertheless, in order to make the desired estimate
rigorous we perform a detailed procedure. Thus, let
$0<\varepsilon<1$ and $\gamma_\varepsilon(\cdot)$ be defined by
\begin{equation} \label{gammaepsilon}
\gamma_\varepsilon(s):=\left\{\begin{aligned}
&1-\frac{1}{s}\quad\text{if }s\geq\varepsilon\\
&1-\frac{1}{\varepsilon}\quad\text{if }s<\varepsilon
\end{aligned}\right.
\end{equation}
Next, we introduce the primitive function $\Gamma_\varepsilon$
defined by
\begin{equation} \label{Gammaepsilon}
\Gamma_\varepsilon(u)=\int_1^u \gamma_\varepsilon (s)ds+1,
\end{equation}
so that it results
\begin{align}
&\Gamma_\varepsilon(u)=u-\log u \quad\text{if }u\geq\varepsilon \notag \\
&\Gamma_\varepsilon(u)=1+\log
\frac{1}{\varepsilon}+\left(1-\frac{1}{\varepsilon}\right)u
\quad\text{if } u<\varepsilon. \label{Gammaepsiloncomputed}
\end{align}
Now, we rewrite \eqref{3approx} in terms of the variables $\Chi^i$
and $U^i$ thus obtaining
\begin{equation}\label{3approxbis}
\frac{U^i-U^{i-1}}{\tau}+\B(U^i(1+\Chi^i))=0.
\end{equation}
Hence, we test \eqref{3approxbis} by $\gamma_\varepsilon(U^i)$ and,
thanks to the convexity of $\Gamma_\varepsilon$, we obtain
\begin{align} \label{1.5}
&\frac{1}{\tau}\int_\Omega
\Gamma_\varepsilon(U^i)-\frac{1}{\tau}\int_\Omega
\Gamma_\varepsilon(U^{i-1})+\int_\Omega \abs{\nabla U^i}^2(1+\Chi^i)\gamma'_\varepsilon(U^i)\notag \\
&+\int_\Omega U^i\nabla
U^i\nabla\Chi^i\gamma'_\varepsilon(U^i)+\int_\Gamma
U^i(1+\Chi^i)\gamma_\varepsilon(U^i)\leq0.
\end{align}
By virtue of \eqref{gammaepsilon}, \eqref{1.5} yields
\begin{align} \label{1.6}
&\frac{1}{\tau}\int_\Omega
\Gamma_\varepsilon(U^i)+\int_{\Omega\cap\{U^i\geq\varepsilon\}}
\abs{\nabla \log U^i}^2(1+\Chi^i)+\int_{\Gamma\cap\{U^i\geq0\}}U^i\notag\\
&\leq\frac{1}{\tau}\int_\Omega
\Gamma_\varepsilon(U^{i-1})+\int_{\Gamma\cap\{U^i\geq0\}}(1+\Chi^i)-\int_{\Omega\cap\{U^i\geq\varepsilon\}}\nabla
\log U^i\nabla\Chi^i.
\end{align}
The second integral on the right hand side turns out to be uniformly
bounded thanks to \eqref{chibounded}, while we can exploit the Young
inequality in order to treat the last term of \eqref{1.6}.
Eventually, using once more \eqref{chibounded}, from \eqref{1.6} we
can obtain
\begin{align} \label{1.7}
&\frac{1}{\tau}\int_\Omega
\Gamma_\varepsilon(U^i)+\frac{1}{2}\int_{\Omega\cap\{U^i\geq\varepsilon\}}
\abs{\nabla \log U^i}^2+\int_{\Gamma\cap\{U^i\geq0\}}U^i\notag\\
&\leq C_0+\frac{1}{\tau}\int_\Omega
\Gamma_\varepsilon(U^{i-1})+\frac{1}{2}\int_{\Omega}\abs{\nabla\Chi^i}^2.
\end{align}
where $C_0$ depends on $\lambda_\beta$. Let us note that, by
inductive hypothesis, we are assuming $\log P^{i-1}\in H$, which
implies $\log U^{i-1}\in H$ (see \eqref{U}), whence $U^{i-1}>0$ a.e.
in $\Omega$. Analogously, since $P^{i-1}\in V$ and
\eqref{chibounded} holds, we have also that $U^{i-1}\in H$. In
particular, we can infer that $\left(U^{i-1}-\log U^{i-1}\right)\in
L^1(\Omega)$. Finally, we observe that the last integral on the
right hand side of \eqref{1.7} is bounded thanks to the regularity
of $\Chi^i$ ($\Chi^i\in W$). Thus, to pass to the limit as
$\varepsilon\searrow0$ in \eqref{1.7}, we can apply the monotone
convergence theorem and get
\begin{align} \label{1.8}
&\frac{1}{\tau}\int_\Omega (U^i-\log U^i)+\frac{1}{2}\int_{\Omega}
\abs{\nabla \log U^i}^2+\int_{\Gamma}U^i\notag\\
&\leq C_0+\frac{1}{\tau}\int_\Omega (U^{i-1}-\log
U^{i-1})+\frac{1}{2}\int_{\Omega}\abs{\nabla\Chi^i}^2.
\end{align}
Due to the previous estimates and to the Poincar\'e-Wirtinger
inequality, from \eqref{1.8} we infer
\begin{equation} \label{tesiinduzionebis}
(U^i-\log U^i)\in L^1(\Omega)\quad \text{and}\quad \log U^i\in V,
\end{equation}
whence $U^i>0$ a.e. in $\Omega$. Now, we combine the estimates
\eqref{1.2} and \eqref{1.8} and sum them. By use of \eqref{U} we
easily deduce
\begin{align} \label{1.9}
&\frac{1}{2\tau}\norm{\Chi^i-\chi_\ast}^2_H+\frac{1}{2\tau}\norm{\Chi^i-\Chi^{i-1}}_H^2-\frac{1}{2\tau}\norm{\Chi^{i-1}-\chi_\ast}^2_H\notag\\
&+\left(\frac{1}{2}+\frac{\nu}{2\tau}\right)\norm{\nabla\Chi^i}^2_H+\frac{\nu}{2\tau}\norm{\nabla\Chi^i-\nabla\Chi^{i-1}}_H^2
-\frac{\nu}{2\tau}\norm{\nabla\Chi^{i-1}}^2_H\notag \\
&+\frac{1}{\tau}\int_\Omega (U^i-\log
U^i)+\frac{1}{2}\int_{\Omega}\abs{\nabla\log U^i}^2+\int_{\Gamma}U^i\notag\\
&+\int_\Omega \log(1+\Chi^{i-1})\Chi^i\leq
C_0+\frac{1}{\tau}\int_\Omega (U^{i-1}-\log
U^{i-1})\notag\\
&+\int_\Omega h (\Theta^{i-1})(\Chi^i-\chi_\ast)-\int_\Omega \log
U^{i-1}(\Chi^i-\chi_\ast)+\int_\Omega \log (1+\Chi^{i-1})\chi_\ast.
\end{align}
Let us observe that the last integral on the left hand side is
non-negative (see \eqref{chibounded}). Besides, owing to
\eqref{chibounded}, we can handle the last integral on the right
hand side of \eqref{1.9} as follows
\begin{equation}\label{ineq2}
\int_\Omega \log (1+\Chi^{i-1}) \chi_\ast\leq\int_\Omega
\Chi^{i-1}\chi_\ast\leq \lambda^2_\beta\abs{\Omega}.
\end{equation}
Hence, by summing up \eqref{1.9} for $i=1,....,m$, with $m\leq
N=T/\tau$, and multiplying by $\tau$, we obtain
\begin{align} \label{1.11}
&\frac{1}{2}\norm{\Chi^m-\chi_\ast}^2_H
+\frac{1}{2}\sum_{i=1}^m\tau\norm{\nabla\Chi^i}^2_H+\frac{\nu}{2}\norm{\nabla\Chi^m}^2_H \notag\\
&+\int_\Omega (U^m-\log U^m)+\frac{1}{2}\sum_{i=1}^m
\tau\int_{\Omega}\abs{\nabla\log
U^i}^2+\sum_{i=1}^m\tau\int_{\Gamma} U^i \notag\\
&\leq\tilde{C_0}+\frac{1}{2}\norm{\Chi^0-\chi_\ast}^2_H+\frac{\nu}{2}\norm{\nabla\Chi^0}^2_H+\int_\Omega(U^0-\log U^0) \notag\\
&+\sum_{i=1}^m \tau\int_\Omega
h(\Theta^{i-1})(\Chi^i-\chi_\ast)-\sum_{i=1}^m\tau\int_\Omega \log
U^{i-1}(\Chi^i-\chi_\ast),
\end{align}
where $\tilde{C_0}=(C_0+\lambda^2_\beta\abs{\Omega})T$. Now, owing
to \eqref{h1} and \eqref{chibounded}, by use of the H\"{o}lder and
Young inequalities we can find two positive constants $C_1$ and
$C_2$ such that the following estimate holds
\begin{align} \label{1.11}
&\frac{1}{2}\norm{\Chi^m-\chi_\ast}^2_H
+\frac{1}{2}\sum_{i=1}^m\tau\norm{\nabla\Chi^i}^2_H+\frac{\nu}{2}\norm{\nabla\Chi^m}^2_H\notag\\
&+\int_\Omega (U^m-\log U^m)+\frac{1}{2}\sum_{i=1}^m
\tau\int_{\Omega}\abs{\nabla\log
U^i}^2+\sum_{i=1}^m\tau\int_{\Gamma}
U^i\notag\\
&\leq
C_1+\frac{1}{2}\norm{\Chi^0-\chi_\ast}^2_H+\frac{\nu}{2}\norm{\nabla\Chi^0}^2_H\notag
\\
&+\int_\Omega U^0+C_2\int_\Omega \abs{\log
U^0}+2\lambda_\beta\sum_{i=1}^{m-1}\tau\norm{\log
U^{i}}_{L^1(\Omega)}.
\end{align}
Let us remind that $U^0=u_{0\tau}=p_{0\tau}/(1+\chi_0)$ and $\log
U^0=\log u_{0\tau}=\log p_{0\tau}-\log(1+\chi_0)$. Then, in virtue
of \eqref{boundchi0} and \eqref{boundp0} we can write
\begin{align} \label{1.12}
&\frac{1}{2}\norm{\Chi^m-\chi_\ast}^2_H
+\frac{1}{2}\sum_{i=1}^m\tau\norm{\nabla\Chi^i}^2_H+\frac{\nu}{2}\norm{\nabla\Chi^m}^2_H\notag\\
&+\int_\Omega(U^m-\log U^m)+\frac{1}{2}\sum_{i=1}^m
\tau\int_{\Omega}\abs{\nabla\log
U^i}^2+\sum_{i=1}^m\tau\int_{\Gamma} U^i\notag\\
&\leq
C_3+\frac{1}{2}\norm{\Chi^0-\chi_\ast}^2_H+\frac{\nu}{2}\norm{\nabla\Chi^0}^2_H
+2\lambda_\beta\sum_{i=1}^{m-1}\tau\norm{\log U^{i}}_{L^1(\Omega)},
\end{align}
for some positive constant $C_3$ not depending on $\tau$. Moreover,
let us note that the following inequality holds
\begin{equation} \label{1.13}
\frac{1}{3}(U^i+\abs{\log U^i})\leq(U^i-\log U^i).
\end{equation}
As a consequence of \eqref{1.13} and owing to \eqref{condiniz} we
can finally obtain
\begin{align} \label{1.14}
&\frac{1}{2}\norm{\Chi^m-\chi_\ast}^2_H
+\frac{1}{2}\sum_{i=1}^m\tau\norm{\nabla\Chi^i}^2_H+\frac{\nu}{2}\norm{\nabla\Chi^m}^2_H\notag\\
&+\frac{1}{3}\norm{U^m}_{L^1(\Omega)} +\frac{1}{3}\norm{\log
U^m}_{L^1(\Omega)}+\frac{1}{2}\sum_{i=1}^m
\tau\int_{\Omega}\abs{\nabla\log
U^i}^2+\sum_{i=1}^m\tau\int_{\Gamma} U^i\notag\\
&\leq
C_3+\frac{1}{2}\norm{\chi_0-\chi_\ast}^2_H+\frac{\nu}{2}\norm{\nabla\chi_0}^2_H
+2\lambda_\beta\sum_{i=1}^{m-1}\tau\norm{\log U^{i}}_{L^1(\Omega)}.
\end{align}
Thus, the Poincar\'e inequality and the discrete Gronwall lemma
\cite[Prop. 2.2.1]{je} applied to \eqref{1.14} imply (see
\eqref{chi0} and \eqref{boundchi0})
\begin{align} \label{estimate1}
&\norm{\bar{\chi}_\tau}_{L^\infty(0,T,V)}+\norm{\bar{u}_\tau}_{L^\infty(0,T,L^1(\Omega))}\notag\\
&+\norm{{\bar{u}_\tau|}_\Gamma}_{L^1(0,T,L^1(\Gamma))}+\norm{\log\bar{u}_\tau}_{L^\infty(0,T,L^1(\Omega))\cap
L^2(0,T,V)} \leq C.
\end{align}
In particular, note that $\bar{u}_\tau>0$ almost everywhere in $Q$.
\\
\\
\textbf{Second a priori estimate.} From the estimates on
$\bar{u}_\tau$ in \eqref{estimate1} it is possible to deduce
ana\-lo\-gous regularities for
$\bar{p}_\tau=\bar{u}_\tau(1+\bar{\chi}_\tau)$. In particular, let
us discuss the regularity of $\log\bar{p}_\tau$. Since
\eqref{chibounded} holds, we can infer that $\bar{\chi}_\tau\geq0$
a.e. in $Q$ and $\bar{\chi}_\tau\in L^\infty(Q)$ (see \eqref{chiinD}
and \eqref{beta^property2}). Hence, thanks to \eqref{defutau}, we
have that
\begin{align} \label{estimate2}
&\norm{\log \bar{p}_\tau}_{L^\infty(0,T,L^1(\Omega))\cap
L^2(0,T,V)}\notag\\
&\leq C_4\left( \norm{\log
\bar{u}_\tau}_{L^\infty(0,T,L^1(\Omega))\cap
L^2(0,T,V)}+\norm{\bar{\chi}_\tau}_{L^\infty(0,T,L^1(\Omega))\cap
L^2(0,T,V)}\right)\leq C.
\end{align}
\\
\\
\textbf{Third a priori estimate.} We test \eqref{2approx} by
$\Chi^i-\Chi^{i-1}$. Thanks to \eqref{eq} and \eqref{2approx'}, we
get
\begin{align} \label{3.1}
&\tau\norm{\frac{\Chi^i-\Chi^{i-1}}{\tau}}^2_H+\nu\tau
\norm{\frac{\nabla(\Chi^i-\Chi^{i-1})}{\tau}}^2_H+\frac{1}{2}\norm{\nabla\Chi^i}^2_H\notag\\
&+\frac{\tau^2}{2}\norm{\frac{\nabla(\Chi^i-\Chi^{i-1})}{\tau}}^2_H-\frac{1}{2}\norm{\nabla\Chi^{i-1}}^2_H
+\int_\Omega \hat{\beta}(\Chi^i)-\int_\Omega\hat{\beta}(\Chi^{i-1})\notag\\
&\leq\int_\Omega (\Chi^i-\Chi^{i-1}) (h(\Theta^{i-1})-\log P^{i-1}).
\end{align}
Now, we can handle the integral on the right hand side of
\eqref{3.1} by means of the Young inequality and get
\begin{align} \label{3.2}
&\int_\Omega(\Chi^i-\Chi^{i-1}) (h(\Theta^{i-1})-\log
P^{i-1})\notag\\
&\leq \frac{\tau}{2}
\norm{\frac{\Chi^i-\Chi^{i-1}}{\tau}}^2_H+\frac{\tau}{2}\left(\abs{\Omega}\norm{h}^2_{L^\infty(\R)}+\norm{\log
P^{i-1}}^2_H\right).
\end{align}
Thus, we add \eqref{3.1} for $i=1,...,m$, with $m\leq N$ and we
easily recover (cf. \eqref{h2}, \eqref{condiniz} and
\eqref{estimate2})
\begin{align} \label{3.4}
&\frac{1}{2}\sum_{i=1}^m\tau\norm{\frac{\Chi^i-\Chi^{i-1}}{\tau}}^2_H+\nu\sum_{i=1}^m\tau
\norm{\frac{\nabla(\Chi^i-\Chi^{i-1})}{\tau}}^2_H+\frac{1}{2}\norm{\nabla\Chi^m}^2_H\notag\\
&+\int_\Omega \hat{\beta}(\Chi^m)\leq
C_5+\frac{1}{2}\norm{\nabla\chi_0}^2_H+\int_\Omega\hat{\beta}(\chi_0).
\end{align}
Hence, owing to \eqref{chi0}, we obtain the following bound
\begin{equation} \label{estimate3}
\norm{\chi_\tau}_{H^1(0,T;V)}+\norm{\hat{\beta}(\bar{\chi}_\tau)}_{L^\infty(0,T;L^1(\Omega))}\leq
C.
\end{equation}
\\
\\
\textbf{Fourth a priori estimate.} We proceed by formally testing
\eqref{2approx} by $\tau\A\Chi^i$. Due to the monotonicity of the
operator $\partial_H J_H$, it turns out that $\tau\int_\Omega\Xi^i\A
\Chi^i\geq0$. Hence, by exploiting once more \eqref{eq} and Young's
inequality, by simi\-lar\-ly proceeding as for \eqref{3.1}, we owe
to \eqref{h2}, \eqref{estimate1} and write
\begin{align} \label{4.1}
&\frac{1}{2}\norm{\nabla\Chi^i}_H^2-\frac{1}{2}\norm{\nabla\Chi^{i-1}}_H^2+\frac{1}{2}\norm{\A\Chi^i}_H^2
-\frac{1}{2}\norm{\A\Chi^{i-1}}_H^2+\tau\norm{\A\Chi^i}_H^2\notag\\
&\leq\frac{\tau}{2}\norm{\A\Chi^i}_H^2+\frac{\tau}{2}\left(\abs{\Omega}\norm{h}^2_{L^\infty(\R)}+\norm{\log
P^{i-1}}^2_H\right).
\end{align}
If we sum up in \eqref{4.1} for $i=1,...,m$, we get (see \eqref{h2},
\eqref{condiniz} and \eqref{estimate2})
\begin{equation} \label{4.2}
\norm{\nabla\Chi^m}_H^2+\norm{\A\Chi^m}_H^2
+\sum_{i=1}^m\tau\norm{\A\Chi^i}_H^2\leq\norm{\nabla\chi_0}_H^2+\norm{\A\chi_0}_H^2+C_6,
\end{equation}
for any $m\leq N$. Finally, the regularity assumptions \eqref{chi0}
on $\chi_0$ ensure that
\begin{equation} \label{estimate4}
\norm{\bar{\chi}_\tau}_{L^\infty(0,T,W)}\leq C.
\end{equation}
\\
\\
\textbf{Fifth a priori estimate.} Test equation \eqref{1approx} by
$E^i$. Ex\-ploiting once more \eqref{eq} and recalling that
$E^i=\psi(\Theta^i,\Chi^i)$ (see \eqref{E}), we can write
\begin{align} \label{5.1}
&\frac{1}{2\tau}\norm{E^i}^2_H+\frac{1}{2\tau}\norm{E^i-E^{i-1}}^2_H-\frac{1}{2\tau}\norm{E^{i-1}}^2_H
\notag\\
&+\int_\Omega
\abs{\nabla\Theta^i}^2\partial_1\psi(\Theta^i,\Chi^i)+\int_\Omega
\nabla\Theta^i\nabla\Chi^i\partial_2\psi(\Theta^i,\Chi^i)\notag\\
&=-\int_\Omega
h(\Theta^{i-1})E^i\frac{\Chi^i-\Chi^{i-1}}{\tau}+\int_\Omega E^i
\left(\frac{\Chi^i-\Chi^{i-1}}{\tau}\right)^2.
\end{align}
To handle the last integral on the right-hand side we use the
extended H\"{o}lder inequality and get
\begin{equation} \label{5.2}
\int_\Omega E^i \left(\frac{\Chi^i-\Chi^{i-1}}{\tau}\right)^2\leq
\norm{E^i}_H \norm{\frac{\Chi^i-\Chi^{i-1}}{\tau}}^2_{L^4(\Omega)}.
\end{equation}
We recall the bounds on $\partial_1\psi$ and $\abs{\partial_2\psi}$
stated in \eqref{derivatapsi1}-\eqref{derivatapsi2}. Then, due to
\eqref{h2}, \eqref{5.2} and the Young and H\"{o}lder inequalities,
we can infer
\begin{align} \label{5.5}
&\frac{1}{2\tau}\norm{E^i}^2_H+\frac{\tau}{2}\norm{\frac{E^i-E^{i-1}}{\tau}}^2_H-\frac{1}{2\tau}\norm{E^{i-1}}^2_H
+\frac{c_s}{2}\norm{\nabla\Theta^i}^2_H \notag\\
&\leq
C_7\norm{\nabla\Chi^i}^2_H+c_h\norm{E^i}_H\norm{\frac{\Chi^i-\Chi^{i-1}}{\tau}}_H+\norm{E^i}_H
\norm{\frac{\Chi^i-\Chi^{i-1}}{\tau}}^2_{L^4(\Omega)},
\end{align}
where $C_7=c_h^2/(2c_s)$. Then, multiplying \eqref{5.5} by $\tau$
and summing up for $i=1,...,m$, we get
\begin{align} \label{5.7}
&\frac{1}{2}\norm{E^m}^2_H+\frac{c_s}{2}\sum_{i=1}^m\tau\norm{\nabla\Theta^i}^2_H\leq\frac{1}{2}\norm{E^0}^2_H+C_7\sum_{i=1}^m\tau\norm{\nabla\Chi^i}^2_H\notag\\
&+\sum_{i=1}^m\norm{E^i}_H\left(c_h\tau\norm{\frac{\Chi^i-\Chi^{i-1}}{\tau}}_H
+\tau\norm{\frac{\Chi^i-\Chi^{i-1}}{\tau}}^2_{L^4(\Omega)}\right).
\end{align}
In order to apply the discrete Gronwall lemma \cite[Prop.
2.2.1]{je}, we have to treat the $m$-th term in the last summation
of \eqref{5.7} separately. Here, we sketch out such a procedure.
Namely, we deal with the $m$-th term by use of the Young inequality
and get
\begin{align} \label{5.8}
&\norm{E^m}_H\left(c_h\tau\norm{\frac{\Chi^m-\Chi^{m-1}}{\tau}}_H
+\tau\norm{\frac{\Chi^m-\Chi^{m-1}}{\tau}}^2_{L^4(\Omega)}\right)\notag \\
&\leq\frac{1}{4}\norm{E^m}^2_H+2\tau
c_h^2\norm{\partial_t\chi_\tau}^2_{L^2(0,T,H)}+2\norm{\partial_t\chi_\tau}^4_{L^2(0,T,L^4(\Omega))}.
\end{align}
Finally, thanks to the Third estimate (see \eqref{estimate3}) and to
the continuous embedding $V\subset L^4(\Omega)$, we conclude that
there exists a positive constant $C_8$ such that
\begin{equation} \label{5.9}
\norm{E^m}_H\left(c_h\tau\norm{\frac{\Chi^m-\Chi^{m-1}}{\tau}}_H
+\tau\norm{\frac{\Chi^m-\Chi^{m-1}}{\tau}}^2_{L^4(\Omega)}\right)\leq
C_8+\frac{1}{4}\norm{E^m}^2_H.
\end{equation}
Owing to \eqref{condinize}, \eqref{estimate1} and \eqref{5.9},
\eqref{5.7} yields
\begin{align} \label{5.10}
&\frac{1}{4}\norm{E^m}^2_H+\frac{c_s}{2}\sum_{i=1}^m\tau\norm{\nabla\Theta^i}^2_H\leq C_9+\frac{1}{2}\norm{e_0}^2_H\notag\\
&+\sum_{i=1}^{m-1}\tau\norm{E^i}_H\left(\frac{c_h^2}{2}+\frac{1}{2}\norm{\frac{\Chi^i-\Chi^{i-1}}{\tau}}^2_H
+\norm{\frac{\Chi^i-\Chi^{i-1}}{\tau}}^2_{L^4(\Omega)}\right).
\end{align}
Eventually, \eqref{e0} and \eqref{estimate1}-\eqref{estimate3} allow
us to make use of the discrete Gronwall lemma \cite[Prop.
2.2.1]{je}, thus obtaining
\begin{equation} \label{estimate5}
\norm{\bar{e}_\tau}_{L^\infty(0,T,H)}+\norm{\nabla\bar{\theta}_\tau}_{L^2(0,T,H)}\leq
C.
\end{equation}
We recall that $\bar{e}_\tau$ is defined as $\bar{e}_\tau=\psi
(\bar{\theta}_\tau, \bar{\chi}_\tau)$ (see \eqref{1tau'}). Thus, the
relation between $\bar{\theta}_\tau$ and $\bar{e}_\tau$ is
bi-lipschitz (cf. \eqref{derivatapsi1}-\eqref{derivatapsi2}) and
allows us to infer from \eqref{estimate5} that
\begin{equation} \label{estimate5bis}
\norm{\bar{\theta}_\tau}_{L^\infty(0,T,H)\cap L^2(0,T,V)}\leq C.
\end{equation}
\\
\\
\textbf{Sixth a priori estimate.} Here, we aim at deriving some
further regularities for the variables $e_\tau$ and $\bar{e}_\tau$
from the previous estimates. At first, let us observe that
\eqref{estimate5bis}, combined with
\eqref{derivatapsi1}-\eqref{derivatapsi2} and with the regularity of
$\bar{\chi}_\tau$ (see \eqref{estimate1}), yields
\begin{equation} \label{6.1}
\norm{\bar{e}_\tau}_{L^2(0,T,V)}\leq C
\end{equation}
Indeed, since $\bar{e}_\tau=\psi
(\bar{\theta}_\tau,\bar{\chi}_\tau),$ the following inequalities
hold
\begin{align}
\norm{\bar{e}_\tau}^2_{L^2(0,T,V)}&=\int_0^T
\norm{\bar{e}_\tau}^2_H+\int_0^T\norm{\nabla\bar{\theta}_\tau\partial_1\psi(\bar{\theta}_\tau,\bar{\chi}_\tau)+\nabla\bar{\chi}_\tau\partial_2\psi(\bar{\theta}_\tau,\bar{\chi}_\tau)}^2_H\notag\\
&\leq
T\norm{\bar{e}_\tau}^2_{L^\infty(0,T,H)}+C_e\norm{\nabla\bar{\theta}_\tau}^2_{L^2(0,T,H)}+c_h^2\norm{\nabla\bar{\chi}_\tau}_{L^2(0,T,H)}\leq
C,
\end{align}
where $C_e$ depends on $c_e$ (see \eqref{derivatapsi1}). Moreover,
from \eqref{6.1} and \eqref{e0} we can also deduce that
\begin{equation} \label{6.2}
\norm{e_\tau}_{L^2(0,T,V)}\leq C.
\end{equation}
Finally, let us observe that \eqref{estimate5bis} yields
\begin{equation}
\norm{\A\bar{\theta}_\tau}_{L^2(0,T,V')}\leq C.
\end{equation}
In addition, from \eqref{estimate3} we can easily deduce that
\begin{equation}
\norm{(\partial_t\chi_\tau)^2}_{L^1(0,T,H)} \quad \mbox{and}\quad
\norm{h(\mathcal{T_\tau}\bar{\theta}_\tau)\partial_t
\chi_\tau}_{L^2(0,T,H)}\leq C.
\end{equation}
Thus, by a comparison in \eqref{1tau}, we can infer that
\begin{equation} \label{6.3}
\norm{\partial_t e_\tau}_{L^1(0,T,V')}\leq C.
\end{equation}
\\
\\
\textbf{Seventh a priori estimate.} We consider \eqref{3approx} and
we rewrite it in terms of $P^i$, $P^{i-1}$ and $\Chi^i$,
$\Chi^{i-1}$, thus obtaining
\begin{equation} \label{3approxnewbis}
\tau^{-1}\frac{P^i}{1+\Chi^i}-\tau^{-1}\frac{P^{i-1}}{1+\Chi^{i-1}}+\B
P^i=0.
\end{equation}
Now, we test \eqref{3approxnewbis} by $P^i$. Thanks to
\eqref{chibounded}, there exists a positive constant
$C_\chi=1/(1+\lambda_\beta)$ such that
\begin{equation} \label{boundchi}
C_\chi\leq\frac{1}{1+\Chi^i}\leq 1 \quad\forall i=0,...,N.
\end{equation}
By standard algebraic calculations, we get
\begin{align} \label{7.1}
&\frac{1}{2\tau} \int_\Omega
\frac{{(P^i)}^2}{1+\Chi^i}-\frac{1}{2\tau} \int_\Omega
\frac{{(P^{i-1})}^2}{1+\Chi^{i-1}}+\frac{1}{2\tau} \int_\Omega
\frac{{(P^i-P^{i-1})}^2}{1+\Chi^i}+\int_\Omega {\abs{\nabla P^i}}^2
+\int_\Gamma (P^i)^2\notag \\
&=\frac{1}{\tau}\int_\Omega P^i
P^{i-1}\frac{\Chi^i-\Chi^{i-1}}{(1+\Chi^i)(1+\Chi^{i-1})}-\frac{1}{2\tau}\int_\Omega
{(P^{i-1})}^2 \frac{\Chi^i-\Chi^{i-1}}{(1+\Chi^i)(1+\Chi^{i-1})}.
\end{align}
Next, we observe that we can find a positive constant $C_P$ such
that
\begin{equation} \label{poincare}
C_P\norm{P^i}^2_V\leq \int_\Omega\abs{\nabla P^i}^2+\int_\Gamma
(P^i)^2.
\end{equation}
Hence, we multiply \eqref{7.1} by $\tau$ and get
\begin{align}\label{7.2}
&\frac{1}{2} \int_\Omega \frac{{(P^i)}^2}{1+\Chi^i}-\frac{1}{2}
\int_\Omega \frac{{(P^{i-1})}^2}{1+\Chi^{i-1}}+\frac{\tau^2}{2}
\int_\Omega
{\left(\frac{P^i-P^{i-1}}{\tau}\right)}^2\frac{1}{1+\Chi^{i}}+C_P\tau\norm{P^i}_V^2\notag\\
&\leq\int_\Omega P^i
P^{i-1}\frac{\Chi^i-\Chi^{i-1}}{(1+\Chi^i)(1+\Chi^{i-1})}-\frac{1}{2}\int_\Omega
{(P^{i-1})}^2 \frac{\Chi^i-\Chi^{i-1}}{(1+\Chi^i)(1+\Chi^{i-1})}.
\end{align}
To handle the first integral on the right-hand side of \eqref{7.2}
we exploit \eqref{boundchi}, the extended H\"{o}lder inequality and
Young's inequality to obtain
\begin{equation} \label{7.3}
\int_\Omega P^i P^{i-1}
\frac{(\Chi^i-\Chi^{i-1})}{(1+\Chi^i)(1+\Chi^{i-1})}\leq \frac{C_P
\tau}{4}\norm{P^i}_V^2+
\frac{\tau}{C_P}\norm{\frac{\Chi^i-\Chi^{i-1}}{\tau}}_V^2
\norm{P^{i-1}}_H^2.
\end{equation}
We proceed similarly for the last integral on the right hand side of
\eqref{7.2} so to get
\begin{equation} \label{7.4}
-\frac{1}{2}\int_\Omega
\frac{{(P^{i-1})}^2(\Chi^i-\Chi^{i-1})}{(1+\Chi^i)(1+\Chi^{i-1})}\leq\frac{C_P\tau}{4}\norm{P^{i-1}}^2_V+\frac{\tau}{4
C_P}\norm{\frac{\Chi^i-\Chi^{i-1}}{\tau}}^2_V\norm{P^{i-1}}^2_H.
\end{equation}
Then, combining \eqref{7.2}-\eqref{7.4} and summing up for
$i=1,...m$ with $m\leq N$, we obtain (see \eqref{condiniz})
\begin{align}\label{7.6}
&\frac{C_\chi}{2} \norm{P^m}_H^2+\frac{C_\chi}{2}\sum_{i=1}^m \tau^2
\norm{\frac{P^i-P^{i-1}}{\tau}}_H^2+\frac{C_P}{2}\sum_{i=1}^m\tau\norm{P^i}_V^2\notag \\
&\leq \frac{1}{2}
\norm{p_{0\tau}}_H^2+\frac{5}{4C_P}\sum_{i=1}^m{\tau\norm{\frac{\Chi^i-\Chi^{i-1}}{\tau}}_V^2
\norm{P^{i-1}}_H^2}.
\end{align}
Hence, we can apply the discrete Gronwall lemma \cite[Prop.
2.2.1]{je} and, due to \eqref{boundp0}, \eqref{estimate3}, deduce
that
\begin{equation} \label{estimate7}
\norm{\bar{p}_\tau}_{L^\infty(0,T,H)\cap L^2(0,T,V)} \leq C.
\end{equation}
\\
\\
\textbf{Eighth a priori estimate.} We test \eqref{3approxnewbis} by
($P^i-P^{i-1}$). By exploiting once more \eqref{eq}, H\"{o}lder's
and Young's inequalities and similarly proceeding as for
\eqref{6.1}, we write
\begin{align} \label{8.1}
&\frac{1}{\tau}\int_\Omega
\frac{(P^i-P^{i-1})^2}{1+\Chi^i}+\frac{1}{2}\norm{\nabla
P^i}^2_H-\frac{1}{2}\norm{\nabla
P^{i-1}}^2_H+\frac{1}{2}\norm{P^i}^2_{L^2(\Gamma)}\notag\\
&-\frac{1}{2}\norm{P^{i-1}}^2_{L^2(\Gamma)}\leq\frac{1}{\tau}\int_\Omega
P^{i-1}(P^i-P^{i-1})\frac{\Chi^i-\Chi^{i-1}}{(1+\Chi^i)(1+\Chi^{i-1})}.
\end{align}
The integral on the right hand side of \eqref{8.1} can be estimated
as follows (see \eqref{boundchi})
\begin{align} \label{8.2}
&\frac{1}{\tau}\int_\Omega
P^{i-1}(P^i-P^{i-1})\frac{\Chi^i-\Chi^{i-1}}{(1+\Chi^i)(1+\Chi^{i-1})}\notag\\
&\leq\frac{C_\chi\tau}{2}\norm{\frac{P^i-P^{i-1}}{\tau}}^2_H+\frac{\tau}{2C_\chi}\norm{P^{i-1}}^2_V\norm{\frac{\Chi^i-\Chi^{i-1}}{\tau}}^2_V.
\end{align}
By combining \eqref{8.1} and \eqref{8.2}, summing up for $i=1,...,m$
and exploiting \eqref{condiniz} along with \eqref{poincare}, we get
\begin{align} \label{8.4}
&\frac{C_\chi}{2}\sum_{i=1}^m
\tau\norm{\frac{P^i-P^{i-1}}{\tau}}_H^2+\frac{C_P}{2}\norm{P^m}_V^2\leq\frac{1}{2}\norm{\nabla
p_{0\tau}}_H^2+\frac{1}{2}\norm{p_{0\tau}}^2_{L^2(\Gamma)}\notag\\
&+\frac{1}{2C_\chi} \sum_{i=1}^m\tau
\norm{\frac{\Chi^i-\Chi^{i-1}}{\tau}}_V^2 \norm{P^{i-1}}_V^2.
\end{align}
By suitably applying the discrete Gronwall lemma \cite[Prop.
2.2.1]{je} it is not difficult to recover the following estimate
(see \eqref{boundp0})
\begin{equation}\label{estimate8}
\norm{p_\tau}_{H^1(0,T,H)}+\norm{\bar{p}_\tau}_{L^\infty(0,T,V)}\leq
C.
\end{equation}
\\
\\
\textbf{Ninth a priori estimate.} Now, we want to achieve the
following estimate
\begin{equation}\label{estimate9}
\norm{u_\tau}_{H^1(0,T,H)}\leq C.
\end{equation}
To this aim, it suffices to prove that
$\norm{\bar{u}_\tau}^2_{L^2(0,T,H)}\leq C \quad \text{and}\quad
\norm{\partial_t u_\tau}^2_{L^2(0,T,H)}\leq C.$ Thus, in order to
check the first estimate, we exploit \eqref{u0}, \eqref{estimate7}
and have
\begin{equation} \label{9.2}
\norm{\bar{u}_\tau}^2_{L^2(0,T,H)}=\norm{\frac{\bar{p}_\tau}{1+\bar{\chi}_\tau}}^2_{L^2(0,T,H)}\leq
\norm{\bar{p}_\tau}^2_{L^2(0,T,H)}\leq C.
\end{equation}
Now, let us deal with the second estimate. By virtue of \eqref{U},
\eqref{boundchi}, \eqref{estimate3} and \eqref{estimate8}, we can
write
\begin{equation} \label{9.3}
\norm{\partial_t u_\tau}^2_{L^2(0,T,H)}\leq 2 \left(\norm{\partial_t
p_\tau}^2_{L^2(0,T,H)}+\norm{\bar{p}_\tau}^2_{L^\infty(0,T,V)}\norm{\partial_t\chi_\tau}^2_{L^2(0,T,V)}\right)\leq
C.
\end{equation}
Thus, combining \eqref{9.2}-\eqref{9.3} we are able to deduce
\eqref{estimate9}. In particular, since $\norm{\partial_t
u_\tau}_{L^2(0,T,H)}\leq C$, a comparison in \eqref{3tau} allows us
to infer that $\norm{\B\bar{p}_\tau}_{L^2(0,T,H)}\leq C$ and
standard elliptic regularity results yield
\begin{equation}
\label{estimate9'} \norm{\bar{p}_\tau}_{L^2(0,T, H^2(\Omega))}\leq
C.
\end{equation}
In the end, the following estimate
\begin{equation} \label{estimate9''}
\norm{\bar{u}_\tau}_{L^\infty(0,T,V)}\leq C
\end{equation}
is a consequence of \eqref{estimate3} and \eqref{estimate8}, since
(see \eqref{defutau})
\begin{align} \label{9.6}
\norm{\bar{u}_\tau}^2_{L^\infty(0,T,V)}&=\sup_{(0,T)}\left(\norm{\frac{\bar{p}_\tau}{1+\bar{\chi}_\tau}}^2_H+\norm{\frac{\nabla\bar{p}_\tau}{1+\bar{\chi}_\tau}-\frac{\bar{p}_\tau\nabla\bar{\chi}_\tau}{(1+\bar{\chi}_\tau)^2}}^2_H\right)\notag\\
&\leq C_{10}
\norm{\bar{p}_\tau}^2_{L^\infty(0,T,V)}\left(1+\norm{\bar{\chi}_\tau}^2_{L^\infty(0,T,W)}\right)\leq
C.
\end{align}
\\
\section{Convergence results}
In this
section we aim to deduce some convergence results that allow us to
pass to the limit as $\tau$ tends to $0$ in
\eqref{1tau}-\eqref{defutau} and \eqref{chitau0},
\eqref{etau0}-\eqref{utau0} so to conclude the proof of Theorem
\ref{teorexistence}. Since we will obtain a great deal of
convergences, we prefer to retrieve them step by step instead of
presenting all the results in just one proposition. At first, we can
combine the previous estimates \eqref{estimate1}, \eqref{estimate2},
\eqref{estimate3}, \eqref{estimate4}, \eqref{estimate5},
\eqref{estimate5bis}, \eqref{6.1}-\eqref{6.3}, \eqref{estimate7},
\eqref{estimate8}, \eqref{estimate9}, \eqref{estimate9'},
\eqref{estimate9''} so to obtain, on account of
\eqref{tau}-\eqref{bartau},
\begin{align}
&\norm{\bar{\chi}_\tau}_{L^\infty(0,T,W)}+\norm{\chi_\tau}_{H^1(0,T,V)}\leq
c,\label{boundforchi} \\
&\norm{\bar{e}_\tau}_{L^\infty(0,T,H)\cap
L^2(0,T,V)}+\norm{e_\tau}_{W^{1,1}(0,T,V')\cap L^2(0,T,V)}\leq
c\label{boundfore}\\
&\norm{\bar{\theta}_\tau}_{L^\infty(0,T,H)\cap L^2(0,T,V)}\leq c \label{boundforteta}\\
&\norm{\bar{u}_\tau}_{L^\infty(0,T,V)}+\norm{u_\tau}_{H^1(0,T,H)}\leq
c \label{boundforu}\\
&\norm{\bar{p}_\tau}_{L^\infty(0,T,V)\cap
L^2(0,T,H^2(\Omega))}+\norm{p_\tau}_{H^1(0,T,H)}\leq
c\label{boundforp}\\
&\norm{\log\bar{p}_\tau}_{L^\infty(0,T,L^1(\Omega))\cap
L^2(0,T,V)}\leq c \label{boundforlogp}
\end{align}
for $c$ not depending on $\tau\in(0,\hat{\tau})$, for a suitable
$\hat{\tau}>0$.

Now, it remains to pass to the limit in \eqref{1tau}-\eqref{defutau}
and \eqref{chitau0}, \eqref{etau0}-\eqref{utau0} as $\tau\searrow
0$. An easy computation yields
\begin{equation} \label{stimavariabilibar}
\norm{\chi_\tau}_{L^\infty(0,T,W)}\leq
\norm{\chi_0}_W+\norm{\bar{\chi}_\tau}_{L^\infty(0,T,W)}\leq c,
\end{equation}
and analogous estimates hold for $\norm{e_\tau}_{L^\infty(0,T,H)}$,
$\norm{p_\tau}_{L^\infty(0,T,V)}$ and
$\norm{u_\tau}_{L^\infty(0,T,V)}$. By virtue of
\eqref{stimavariabilibar}, well-known weak and weak star compactness
results apply to \eqref{boundforchi}-\eqref{boundfore},
\eqref{boundforu}-\eqref{boundforp} and ensure the following weak
and weak star convergences to hold, possibly for a subsequence of
$\tau$,
\begin{align}
&\chi_\tau \stackrel{\ast}{\rightharpoonup}\chi\quad\text{in }
H^1(0,T,V)\cap L^\infty(0,T,W),\label{convergencechitau}\\
&e_\tau \stackrel{\ast}{\rightharpoonup} e\quad\text{in }
L^\infty(0,T,H)\cap L^2(0,T,V), \label{convergenceetau}\\
&u_\tau\stackrel{\ast}{\rightharpoonup} u \quad\text{in
}H^1(0,T,H)\cap L^\infty(0,T,V),
\label{convergenceutau}\\
&p_\tau \stackrel{\ast}{\rightharpoonup}p\quad\text{in }
H^1(0,T,H)\cap L^\infty(0,T,V).\label{convergenceptau}
\end{align}
Let us stress that, even if we do not specify it, the convergence
results have to be intended to hold up to the extraction of a
suitable subsequence of $\tau$, still denoted by $\tau$ for the sake
of convenience. Then, owing to strong compactness theorems (see e.g.
\cite[Cor. 4, p. 85]{si}), by \eqref{6.3},
\eqref{convergencechitau}-\eqref{convergenceptau}, we get
\begin{align}
&\chi_\tau\rightarrow\chi\quad\text{in }
C^0([0,T],H^{2-\varepsilon}(\Omega)),\label{strongconvergencechitau}\\
&e_\tau\rightarrow e\quad\text{in }L^2(0,T,H)
\label{strongconvergenceetau}\\
&u_\tau \rightarrow u \quad \text{and}\quad p_\tau\rightarrow p
\quad\text{in } C^0([0,T],H^{1-\varepsilon}(\Omega)),\quad\text{if
}\varepsilon>0. \label{strongconvergenceptau}
\end{align}

Moreover, the following relation is fulfilled (cf.
\eqref{tau}-\eqref{bartau})
\begin{equation} \label{differenceLinfty}
\norm{\chi_\tau-\bar{\chi}_\tau}^2_{L^\infty(0,T,V)}\leq\max_{1\leq
i\leq N} \tau^2\norm{\frac{\Chi^i-\Chi^{i-1}}{\tau}}^2_H\leq
\tau\norm{\partial_t \chi_\tau}^2_{L^2(0,T,H)}\leq \tau c,
\end{equation}
and analogous estimates holds for
$\norm{p_\tau-\bar{p}_\tau}^2_{L^\infty(0,T,H)}$ and
$\norm{u_\tau-\bar{u}_\tau}^2_{L^\infty(0,T,H)}$. While for the
difference between $e_\tau$ and $\bar{e}_\tau$ we have
\begin{equation} \label{differenceL1}
\norm{e_\tau-\bar{e}_\tau}_{L^1(0,T,V')}\leq\tau\norm{\partial_t
e_\tau}_{L^1(0,T,V')}\leq \tau c.
\end{equation}
Finally, with the help of \eqref{boundforchi}-\eqref{boundforlogp},
\eqref{strongconvergencechitau}-\eqref{strongconvergenceptau} and
\eqref{differenceLinfty}-\eqref{differenceL1}, we are allowed to
infer that
\begin{align}
&\bar{\chi}_\tau\rightarrow\chi\quad\text{in }L^\infty(0,T,V),
\quad \bar{\chi}_\tau\stackrel{\ast}{\rightharpoonup}\chi\quad\text{in }L^\infty(0,T,W),\label{convergencechibartau}\\
&\bar{e}_\tau\rightarrow e\quad\text{in } L^1(0,T,V'),\quad
\bar{e}_\tau \stackrel{\ast}{\rightharpoonup} e\quad\text{in
}L^\infty(0,T,H)\cap L^2(0,T,V),\label{convergenceebartau}\\
&\bar{\theta}_\tau \stackrel{\ast}{\rightharpoonup} \theta\quad\text{in }L^\infty(0,T,H)\cap L^2(0,T,V),\label{convergencetetabartau}\\
&\bar{u}_\tau\rightarrow u \quad\text{in } L^\infty(0,T,H),\quad
\bar{u}_\tau\stackrel{\ast}{\rightharpoonup}u\quad\text{in
}L^\infty(0,T,V),\label{convergenceubartau}\\
&\bar{p}_\tau\rightarrow p \quad\text{in } L^\infty(0,T,H),\quad
\bar{p}_\tau\stackrel{\ast}{\rightharpoonup}p\quad\text{in
}L^\infty(0,T,V)\cap L^2(0,T,
H^2(\Omega)),\label{convergencepbartau}\\
&\log\bar{p}_\tau\rightharpoonup y \quad\text{in }L^2(0,T,V),
\label{convergencelopbartau1}
\end{align}
where the last limit $y$ will be identified in the sequel.

Now, we want to improve the strong convergence for $\bar{e}_\tau$ in
\eqref{convergenceebartau}. To this aim, we first note that
\begin{equation}\label{stimaperebartauinLpV'}
\norm{\bar{e}_\tau-e}^p_{L^p(0,T,V')}\leq
\norm{\bar{e}_\tau-e}^{p-1}_{L^\infty(0,T,V')}\norm{\bar{e}_\tau-e}_{L^1(0,T,V')},
\end{equation}
where, by virtue of \eqref{convergenceebartau},
$\norm{\bar{e}_\tau-e}^{p-1}_{L^\infty(0,T,V')}\leq c$ and
$\norm{\bar{e}_\tau-e}_{L^1(0,T,V')}\rightarrow 0$ for all $1\leq
p<+\infty$. Thus, we have $\bar{e}_\tau\rightarrow e\quad\text{in }
L^p(0,T,V').$ In addition, we can perform the following estimate
\begin{equation} \label{stimaperebartauinL2H}
\norm{\bar{e}_\tau-e}^2_{L^2(0,T,H)}\leq\norm{\bar{e}_\tau-e}_{L^2(0,T,V)}\norm{\bar{e}_\tau-e}_{L^2(0,T,V')}
\end{equation}
and conclude that
\begin{equation} \label{strongconvergenceebartauinL2H}
\bar{e}_\tau\rightarrow e\quad\text{in } L^2(0,T,H).
\end{equation}
Finally, since the relation between $\bar{\theta}_\tau$ and
$\bar{e}_\tau$ is bi-lipschitz continuous (see
\eqref{derivatapsi1}), from \eqref{strongconvergenceebartauinL2H}
and \eqref{convergencechibartau} we can deduce that
$\bar{\theta}_\tau$ is a Cauchy sequence in $L^2(0,T,H)$ which is a
complete space. Therefore, recalling the convergence
$\bar{\theta}_\tau\stackrel{\ast}{\rightharpoonup}\theta$ in
$L^\infty(0,T,H)\cap L^2(0,T,V)$, by uniqueness of the limit we
conclude that
\begin{equation} \label{strongconvergencetetabartauinL2H}
\bar{\theta}_\tau \rightarrow\theta\quad\text{in } L^2(0,T,H).
\end{equation}
Now, by virtue of the above convergences
\eqref{strongconvergenceebartauinL2H},
\eqref{strongconvergencetetabartauinL2H},
\eqref{convergencechibartau} and owing to the properties of $\psi$
(see \eqref{derivatapsi1}-\eqref{derivatapsi2}), it is easy to
recover \eqref{1'} from \eqref{1tau'}, i.e. $e=\psi(\theta,\chi)$.

Next, we deal with the logarithmic term in \eqref{2tau}. Thanks to
the strong convergence in \eqref{convergencepbartau}, there exists a
subsequence still denoted by $\bar{p}_\tau$ such that
\begin{equation} \label{convergenceqopbartau}
\bar{p}_\tau\rightarrow p \quad \text{a.e. in }Q,
\end{equation}
and consequently such that
\begin{equation} \label{convergenceqologpbartau}
\log\bar{p}_\tau\rightarrow\log p \quad \text{a.e. in }Q.
\end{equation}
In principle, the a.e. limit of $\log \bar{p}_\tau$ could be
$-\infty$ in a subset of positive measure (in which $p=0$): but the
property
$$\int_Q \abs{\log \bar{p}_\tau}^2\leq c\quad\forall\tau$$
and the Fatou lemma imply that
$$\int_Q \abs{\log p}^2\leq \liminf_{\tau\searrow0} \int_Q\abs{ \log \bar{p}_\tau}^2\leq c,$$
whence $\log p$ is well defined. Actually, we are now about to show
that the weak limit $y$ of the sequence $\log \bar{p}_\tau$
coincides with the a.e. limit, i.e. with $\log p$. By virtue of
\eqref{convergenceqologpbartau}, we can invoke the Severini-Egorov
theorem and deduce that for all $\varepsilon>0$ there exists a set
$Q_\varepsilon\subset Q$ such that
$\text{meas}(Q_\varepsilon)<\varepsilon$ and
\begin{equation} \label{convergenceaulogpbartau}
\log\bar{p}_\tau\rightarrow\log p \quad \text{uniformly in
}Q\setminus Q_\varepsilon.
\end{equation}
In addition, since \eqref{boundforlogp} holds, a standard
interpolation calculus and the con\-ti\-nuous embedding $V\subset
L^6(\Omega)$ yield
\begin{align*}
\norm{\log\bar{p}_\tau}^{8/3}_{L^{8/3}(Q)}&\leq
\int_0^T\norm{\abs{\log
\bar{p}_\tau}^{\frac{2}{3}}}_{L^{3/2}(\Omega)}\norm{\abs{\log
\bar{p_\tau}}^2}_{L^3(\Omega)}\notag\\
&\leq\norm{\log\bar{p}_\tau}^{2/3}_{L^\infty(0,T,L^1(\Omega))}\norm{\log\bar{p}_\tau}_{L^2(0,T,L^6(\Omega))}\notag\\
&\leq\norm{\log\bar{p}_\tau}^{2/3}_{L^\infty(0,T,L^1(\Omega))}\norm{\log\bar{p}_\tau}_{L^2(0,T,V)}\leq
c,
\end{align*}
whence
\begin{equation}\label{boundforlogpinQ}
\norm{\log\bar{p}_\tau}_{L^{8/3}(Q)}\leq c.
\end{equation}
By virtue of \eqref{convergenceaulogpbartau}-\eqref{boundforlogpinQ}
and of the H\"{o}lder inequality, if $1\leq q<8/3$ we can write
\begin{align}
&\int_Q\abs{\log\bar{p}_\tau-\log p}^q\leq \left(\int_{
Q_\varepsilon}\abs{\log\bar{p}_\tau-\log
p}^{8/3}\right)^{\frac{3q}{8}}\varepsilon^{1-\frac{3q}{8}}+\int_{Q\setminus
Q_\varepsilon} \abs{\log\bar{p}_\tau-\log p}^q, \notag\\
\intertext{and then}
&\lim_{\tau\searrow0}\int_Q\abs{\log\bar{p}_\tau-\log
p}^q\leq\sup_{\tau} \left(\int_{Q}\abs{\log\bar{p}_\tau-\log
p}^{8/3}\right)^{\frac{3q}{8}}\varepsilon^{1-\frac{3q}{8}}\leq C
\varepsilon^{1-\frac{3q}{8}}.\notag
\end{align}
Taking the limit as $\varepsilon\searrow0$ we obtain
$\log\bar{p}_\tau\rightarrow\log p \quad \text{in
}L^q(Q)\quad\forall 1\leq q<\frac{8}{3}$, and in particular
\begin{equation} \label{finalconvergencelogpbartau}
\log\bar{p}_\tau\rightarrow\log p \quad \text{in }L^2(Q),
\end{equation}
as $\tau\searrow0$. Then, owing to Proposition
\ref{passaggioallimite}, the same convergence easily holds (see
\eqref{bartau} and \eqref{translation}) for
$\mathcal{T_\tau}\log(\bar{p}_\tau)$. Besides, we can identify the
weak limit $y$ (see \eqref{convergencelopbartau1}): the convergence
in \eqref{finalconvergencelogpbartau} definitely yields $y=\log p$.

Now, by the above convergences, we are in the position of taking the
limit as $\tau\searrow0$ in \eqref{1tau}-\eqref{defutau} and
\eqref{chitau0}, \eqref{etau0}-\eqref{utau0}. At first, we observe
that (see \eqref{p0} and \eqref{p0tau})
\begin{equation}\label{convergencep0tau}
p_\tau(0)=p_{0\tau}\rightarrow p_0\quad\text{in }H.
\end{equation}
Thanks to \eqref{chitau0}, \eqref{etau0}-\eqref{utau0},
\eqref{convergencep0tau}, \eqref{convergencechitau} and
\eqref{convergenceutau}, the limit functions $\chi$, $e$, $u$
satisfy the initial conditions \eqref{chiinitial},
\eqref{einitial}-\eqref{uinitial}. The rest of the proof will
proceed in five steps. As first, we will pass to the limit in
\eqref{2tau}; secondly we will prove that the following strong
convergence holds
\begin{equation} \label{strongconvergencechitaut}
\partial_t\chi_\tau\rightarrow\chi_t \quad\text{in }
L^2(0,T,V).
\end{equation}
Then, by virtue of \eqref{strongconvergencechitaut} we will discuss
the passage to the limit in \eqref{1tau}. Finally we will consider
\eqref{3tau} and we will easily recover \eqref{3} by virtue of the
above convergences.

In order to pass to the limit in \eqref{2tau}, we first observe that
\eqref{strongconvergencetetabartauinL2H} combined with \eqref{h1}
yields (see Proposition \ref{passaggioallimite})
\begin{equation} \label{convergencehteta}
h(\mathcal{T_\tau}\bar{\theta}_\tau)\rightarrow
h(\theta)\quad\text{in } L^2(Q).
\end{equation}
Thus, by virtue of \eqref{convergencechitau},
\eqref{convergencechibartau}, \eqref{finalconvergencelogpbartau} and
\eqref{convergencehteta}, in order to conclude the passage to the
limit as $\tau\searrow0$ in \eqref{2tau}, it suffices to control the
sequence $\bar{\xi}_\tau$. In particular we need to verify that
$\bar{\xi}_\tau$ converges, in a suitable sense, to some selection
$\xi\in\partial_{V,V'}J(\chi)$. To this aim, we observe that, by a
comparison in \eqref{2tau}, \eqref{h1}, \eqref{boundforchi} and
\eqref{boundforlogp} imply
\begin{equation} \label{boundforxi}
\norm{\bar{\xi}_\tau}_{L^2(0,T,V')}\leq c,
\end{equation}
and consequently we have
\begin{equation} \label{convergencexibartau}
\bar{\xi}_\tau\rightharpoonup\xi\quad\text{in }L^2(0,T,V').
\end{equation}
Hence, by \eqref{convergencechibartau} and
\eqref{convergencexibartau} we can deduce
\begin{equation} \label{convergenceduality}
\int_0^T\langle\bar{\xi}_\tau,\bar{\chi}_\tau\rangle\rightarrow\int_0^T\langle\xi,\chi\rangle,
\end{equation}
as $\tau\searrow0$, which enables us to apply the result presented
in \cite[Prop. 2.5, p. 27]{bre} for $X=L^2(0,T,V)$ and deduce
\begin{equation} \label{limitinclusion}
\xi\in\partial_{V,V'}J(\chi)\quad\text{a.e. in }(0,T).
\end{equation}
Indeed, it is known that a maximal monotone operator from $V$ to
$V'$ induces an analogous operator from $L^2(0,T,V)$ to
$L^2(0,T,V')$ which is defined by the a.e. relation in $(0,T)$.
Moreover, in our framework we can refer to \cite[Ex. 2.3.3, p.
25]{bre} and deduce that the induced operator from $L^2(0,T,V)$ to
$L^2(0,T,V')$ is maximal monotone as well. Finally, by the above
arguments we can pass to the limit in \eqref{2tau} and get \eqref{2}
solved by $\chi$, $\theta$, $p$ in $V'$ and a.e. in (0,T). Now, we
aim to prove \eqref{strongconvergencechitaut}. We first note that,
since \eqref{convergencechitau} holds we can infer that
\begin{equation} \label{convergencechitaut}
\partial_t\chi_\tau\rightharpoonup\chi_t\quad\text{in }L^2(0,T,V).
\end{equation}
As a consequence, \eqref{strongconvergencechitaut} can be obtained
just by verifying that
\begin{equation} \label{limsup}
\limsup_{\tau\searrow0}
\norm{\partial_t\chi_\tau}^2_{L^2(0,T,V)}\leq\norm{\chi_t}^2_{L^2(0,T,V)},
\end{equation}
since the strong convergence of the norms combined with the weak
convergence imply the required strong convergence
\eqref{strongconvergencechitaut}. To obtain \eqref{limsup} we test
\eqref{2tau} by $\partial_t\chi_\tau$, integrate over $(0,T)$, and
take the $\limsup$ as $\tau\searrow0$. We have
\begin{align} \label{test}
\limsup_{\tau\searrow0}
\int_0^T\norm{\partial_t\chi_\tau}^2_H+\nu\norm{\nabla\partial_t\chi_\tau}^2_H&=\limsup_{\tau\searrow0}\left(
-\int_Q\nabla\bar{\chi}_\tau\nabla\partial_t\chi_\tau-\int_0^T\langle\bar{\xi}_\tau,\partial_t\chi_\tau\rangle\right. \notag\\
&\left.+\int_Q
h(\mathcal{T_\tau}\bar{\theta}_\tau)\partial_t\chi_\tau-\int_Q
\mathcal{T_\tau}(\log\bar{p}_\tau)\partial_t\chi_\tau\right).
\end{align}
We first note that, since \eqref{convergencechitau},
\eqref{convergencechibartau}, \eqref{convergencehteta} and
\eqref{finalconvergencelogpbartau} hold, we can infer that
\begin{align} \label{lim}
&\lim_{\tau\searrow0}\left(-\int_Q\nabla\bar{\chi}_\tau\nabla\partial_t\chi_\tau+\int_Q
h(\mathcal{T_\tau}\bar{\theta}_\tau)\partial_t\chi_\tau-\int_Q
\mathcal{T_\tau}(\log\bar{p}_\tau)\partial_t\chi_\tau\right)\notag\\
&=-\int_Q\nabla\chi\nabla\chi_t+\int_Q h(\theta)\chi_t-\int_Q(\log
p) \chi_t.
\end{align}
Next, we have to treat the term
$-\int_0^T\langle\bar{\xi}_\tau,\partial_t\chi_\tau\rangle$. To this
aim, we observe that by definition of subdifferential and owing to
\eqref{tau}-\eqref{bartau}, we can write
\begin{equation}\label{difficultterm}
\int_0^T \langle\bar{\xi}_\tau,\partial_t\chi_\tau\rangle
=\sum_{i=1}^N\langle\Xi^i,\Chi^i-\Chi^{i-1}\rangle\geq \sum_{i=1}^N
J(\Chi^i)-J(\Chi^{i-1})=J(\chi_\tau(T))-J(\chi_0)
\end{equation}
for any $\tau\geq0$. Hence, by the lower semicontinuity in $V$ of
the function $J$ and due to \eqref{strongconvergencechitau}, we
claim that
\begin{equation} \label{limdifficultterm}
\limsup_{\tau\searrow0}-\int_0^T
\langle\bar{\xi}_\tau,\partial_t\chi_\tau\rangle\leq-J(\chi(T))+J(\chi_0).
\end{equation}
Hence, by combining \eqref{lim} and \eqref{limdifficultterm} we get
\begin{align}\label{limsupfinale}
&\limsup_{\tau\searrow0}\left(-\int_Q\nabla\bar{\chi}_\tau\nabla\partial_t\chi_\tau-\int_0^T\langle\bar{\xi}_\tau,\partial_t\chi_\tau\rangle\right.\left.+\int_Q
h(\mathcal{T_\tau}\bar{\theta}_\tau)\partial_t\chi_\tau-\int_Q
\mathcal{T_\tau}(\log\bar{p}_\tau)\partial_t\chi_\tau\right)\notag\\
&\leq-\int_Q\nabla\chi\nabla\chi_t+\int_Q
h(\theta)\chi_t-\int_Q(\log p)\chi_t-J(\chi(T))+J(\chi_0),
\end{align}
and the right hand side of \eqref{limsupfinale} is equal to
\begin{equation}
\int_0^T\norm{\chi_t}^2_H+\nu\norm{\nabla\chi_t}^2_H,
\end{equation}
as one can easily verify by testing \eqref{2} by $\chi_t$ and then
integrating in time. Thus \eqref{strongconvergencechitaut} is
proved.
\begin{obs}
The last result follows once proved that
\begin{equation}
-J(\chi(T))+J(\chi_0)=-\int_0^T\langle\xi,\chi_t\rangle.
\end{equation}
This can be obtained by extending the statement in \cite[Lemma 3.3,
p. 73]{bre} to the case of abstract subdifferential operators
defined in the duality pairing between $V'$ and $V$.
\end{obs}
Now, we discuss how to perform the passage to the limit in
\eqref{1tau}. By virtue of \eqref{convergencetetabartau} we have
\begin{equation} \label{convergenceAteta}
\A\bar{\theta}_\tau\rightharpoonup\A\theta \quad\text{in }
L^2(0,T,V').
\end{equation}
In addition, owing to \eqref{convergencehteta} and
\eqref{strongconvergencechitaut} we may infer that
\begin{align} \label{convergencetwoterms}
&h(\mathcal{T_\tau}\bar{\theta}_\tau)\partial_t\chi_\tau\rightarrow
h(\theta)\chi_t\quad\text{in
}L^1(0,T,L^{3/2}(\Omega)),\notag\\
&\partial_t\chi_\tau^2\rightarrow\chi_t^2\quad\text{in }L^1(0,T,H).
\end{align}
Thus, by a comparison in \eqref{1tau} we can deduce that
\begin{equation} \label{convergenceetaut}
\partial_te_\tau\rightharpoonup\eta\quad\text{in }L^1(0,T,V'),
\end{equation}
where $\eta=-\A\theta-h(\theta)\chi_t+\chi_t^2$. In order to show
that $\eta=e_t$ we argue as follows. We consider the time
convolution product $1\ast(\partial_t e_\tau)$ which satisfies
\begin{equation} \label{convolutionproductetaut}
1\ast(\partial_t e_\tau)=e_\tau-e_0.
\end{equation}
Then, we observe that the following convergences hold thanks to
\eqref{convergenceetaut}
$$1\ast(\partial_t e_\tau)\rightharpoonup 1\ast \eta\quad\text{in } W^{1,1}(0,T,V'),$$
while
$$e_\tau-e_0\stackrel{\ast}{\rightharpoonup} e-e_0\quad\text{in }L^\infty(0,T,H).$$ The uniqueness of the limit of \eqref{convolutionproductetaut} entails
\begin{equation} \label{limitconvolutionproduct}
e=e_0+1\ast\eta.
\end{equation}
Hence, $e$ must be derivable and, by deriving
\eqref{limitconvolutionproduct}, we get $e_t=\eta$. Thus, we can
pass to the limit in \eqref{1tau} getting \eqref{1}.

Finally, we aim to take the limit in \eqref{3tau}-\eqref{defutau}.
On account of \eqref{convergencepbartau} it is easily verified that
\begin{equation} \label{convergenceBpbartau}
\B\bar{p}_\tau \stackrel{\ast}{\rightharpoonup}\B p\quad\text{in
}L^\infty(0,T,V')\cap L^2(0,T,H).
\end{equation}
On the other hand, thanks to \eqref{convergenceutau}, we have
\begin{equation} \label{convergenceutaut}
\partial_t u_\tau\rightharpoonup u_t\quad\text{in }L^2(0,T,H).
\end{equation}
As a consequence of \eqref{convergenceBpbartau} and
\eqref{convergenceutaut} we can pass to the limit in \eqref{3tau}
and get \eqref{3}. Finally, we equivalently rewrite \eqref{defutau}
as follows
\begin{equation}\label{differenceu}
\bar{u}_\tau=\frac{1}{1+\bar{\chi}_\tau}\left(\bar{p}_\tau-p\right)+
p\left(\frac{1}{1+\bar{\chi}_\tau}\right).
\end{equation}
Let us note that \eqref{convergencechibartau},
\eqref{convergencepbartau} and the fact that $\bar{\chi}_\tau\in
L^\infty(Q)$ allow us to take the limit as $\tau\searrow0$ in
\eqref{differenceu} and eventually get \eqref{defu}. Then Theorem
\ref{teorexistence} is completely proved.

\section{Regularity results}
This section is devoted to the proof of Theorems
\ref{teorregularitypressure}, \ref{teorfurtherregularities} and
\ref{teorpositivitytemperature}. The proof of the improved
regularities is based on formal estimates performed on the solutions
of the complete dissipative model \eqref{1}-\eqref{defu}. Actually,
the following estimates can be made rigorous and, in particular, the
proof of Theorem \ref{teorfurtherregularities} can be directly
reproduced on the discrete scheme.

As first, we prove Theorem \ref{teorregularitypressure} that states
a further regularity for the inverse of the pressure $p^{-1}$. We
first achieve an improved regularity for the variable $u^{-1}$ by
means of a formal estimate, then we will deduce the required result
for $p^{-1}$.

Before proceeding, for $n\leq 3$, we recall the Sobolev embedding
$H^1(\Omega)\hookrightarrow L^6(\Omega)$ and the Gagliardo-Niremberg
inequality (cf. \cite{ni}), yielding, for $n=3$,
\begin{equation} \label{gagliardoniremberg}
\norm{v}^2_{L^3(\Omega)}\leq C_{\text{GN}} \norm{v}_H \norm{\nabla
v}_H+ C'_{\text{GN}} \norm{v}^2_H.
\end{equation}
Now, we consider \eqref{3} rewritten in the terms of the variables
$u$ and $\chi$ as
\begin{equation} \label{3rewritten}
u_t+\B(u(1+\chi))=0.
\end{equation}
We formally test \eqref{3rewritten} by $-u^{-3}$ and integrate over
$(0,t)$. For a rigorous estimate, we could proceed as in the First a
priori estimate of Chapter 4. In particular, we could truncate the
function $-u^{-3}$ at level $\varepsilon$ and then let $\varepsilon$
tend to zero. Let us note that, by virtue of \eqref{uinitial},
\eqref{chi0} and \eqref{p^-10} we get
$u_0^{-1}=(1+\chi_0)p_0^{-1}\in H$. We first have
\begin{equation} \label{p.1}
-\int_0^t\int_\Omega u_t u^{-3}=\frac{1}{2}\int_0^t\int_\Omega
\frac{\partial}{\partial
t}(u^{-2})=\frac{1}{2}\norm{u^{-1}(t)}^2_H-\frac{1}{2}\norm{u_0^{-1}}^2_H.
\end{equation}
Then, by definition of $\B$, we write
\begin{align} \label{p.2}
&-\int_0^t \langle\B (u(1+\chi)),u^{-3}\rangle
\notag\\
&=3\int_0^t\int_\Omega (1+\chi)\abs{\nabla u}^2
u^{-4}-\int_0^t\int_\Gamma u^{-2} (1+\chi)+3\int_0^t\int_\Omega
u^{-3} \nabla u \nabla \chi.
\end{align}
The third integral on the right hand side of \eqref{p.2} is
estimated as follows (see \eqref{gagliardoniremberg})
\begin{align} \label{p.3}
&3\int_0^t\int_\Omega \abs{u^{-3}\nabla u\nabla\chi}\leq3\int_0^t
\norm{\nabla\chi}_{L^6(\Omega)}\norm{\nabla(u^{-1})}_H\norm{u^{-1}}_{L^3(\Omega)}\notag\\
&\leq \frac{3}{4} \norm{\nabla(u^{-1})}^2_{L^2(0,t,H)}+3
C_{\text{GN}}\int_0^t\norm{\nabla\chi}^2_{L^6(\Omega)}\norm{u^{-1}}_H\norm{\nabla(u^{-1})}_H\notag\\
&+3 C'_{\text{GN}}\int_0^t\norm
{\nabla\chi}^2_{L^6(\Omega)}\norm{u^{-1}}^2_H.
\end{align}
Now, since $\nabla\chi\in L^\infty(0,T,V)$ (cf.
\eqref{chiregularity}), with the help of Young's inequality, we can
eventually write
\begin{equation} \label{p.4}
3\int_0^t\int_\Omega \abs{u^{-3}\nabla u\nabla\chi}\leq
C_{11}+\frac{3}{2}\norm{\nabla(u^{-1})}^2_{L^2(0,t,H)}+C_{12}\int_0^t\norm{u^{-1}}^2_H.
\end{equation}
As the trace operator $\gamma:V\rightarrow L^2(\Gamma)$ is compact,
we may deduce that, for any $\sigma>0$, there exists $c_\sigma>0$
such that
\begin{equation} \label{p.5}
\norm{\gamma(v)}^2_{L^2(\Gamma)}\leq \sigma \norm{v}^2_V+C_\sigma
\norm{v}^2_H \quad \forall v\in V.
\end{equation}
Thus, we can control the boundary integral in \eqref{p.2} as follows
\begin{equation} \label{p.6}
\abs{\int_0^t\int_\Gamma u^{-2} (1+\chi)}\leq \sigma \norm{\nabla
u^{-1}}^2_{L^2(0,t,H)}+C_\sigma \norm{u^{-1}}^2_{L^2(0,t,H)}.
\end{equation}
Finally we point out that for the first integral on the right hand
side of \eqref{p.2} the following equality holds
\begin{equation} \label{p.7}
3\int_0^t\int_\Omega (1+\chi) \abs{\nabla u}^2
u^{-4}=3\int_0^t\int_\Omega (1+\chi)\abs{\nabla (u^{-1})}^2.
\end{equation}
Combining \eqref{3rewritten}-\eqref{p.6} for a sufficiently small
$\sigma$, we have
\begin{equation} \label{p.8}
\norm{u^{-1}(t)}^2_H+\norm{\nabla (u^{-1})}^2_{L^2(0,t,H)}\leq
C_{13}\left(1+\int_0^t \norm{u^{-1}}^2_H\right).
\end{equation}
We can apply the Gronwall lemma \cite[Theorem 2.1]{ba} to
\eqref{p.8} and finally obtain
\begin{equation} \label{furtherubound}
\norm{u^{-1}}_{L^\infty(0,T,H)\cap L^2(0,T,V)}\leq c.
\end{equation}
The same estimate can be derived for the inverse of the pressure
$p^{-1}$. Indeed, thanks to \eqref{chiregularity} and to the
relation $p^{-1}=u^{-1}(1+\chi)^{-1}$, we can easily deduce that
\eqref{furtherubound} implies
\begin{equation} \label{furtherpbound}
\norm{p^{-1}}_{L^\infty(0,T,H)\cap L^2(0,T,V)}\leq c.
\end{equation}
Let us detail such a procedure. By virtue of the H\"{o}lder
inequality and of the continuous embedding $V\subset L^4(\Omega)$,
we have
\begin{align} \label{counts1}
&\norm{p^{-1}}^2_{L^2(0,T,V)}=\int_0^T\norm{\frac{u^{-1}}{1+\chi}}^2_H+\int_0^T\norm{\frac{\nabla
u^{-1}}{1+\chi}-\frac{u^{-1}\nabla\chi }{(1+\chi)^2}}^2_H\notag\\
&\leq C_{14}\left(\norm{u^{-1}}^2_{L^2(0,T,V)}+\int_0^T
\norm{u^{-1}}^2_{L^4(\Omega)}\norm{\nabla\chi}^2_{L^4(\Omega)}\right)\notag\\
&\leq
C_{15}\left(\norm{u^{-1}}^2_{L^2(0,T,V)}+\norm{\chi}^2_{L^\infty(0,T,W)}\norm{u^{-1}}^2_{L^2(0,T,V)}\right)\leq
c.
\end{align}
Analogously, we can easily deduce that
$\norm{p^{-1}}^2_{L^\infty(0,T,H)}\leq
\norm{u^{-1}}^2_{L^\infty(0,T,H)}.$ From \eqref{counts1} we can
eventually infer that \eqref{furtherpbound} holds and Theorem
\ref{teorregularitypressure} is consequently proved.

In order to prove Theorem \ref{teorfurtherregularities} we perform a
further a priori estimate, that can be suitably reproduced on the
discrete scheme. Let us deal with \eqref{1} and rewrite it as
\begin{equation} \label{t.0}
(1+\chi\theta h''(\theta))\theta_t+\A\theta=-\theta
h'(\theta)\chi_t+\chi^2_t.
\end{equation}
Now, we test \eqref{t.0} by $\theta_t$, integrate over $(0,t)$ and
get
\begin{equation}
\int_0^t\int_\Omega \theta_t^2 (1+\chi\theta
h''(\theta))+\int_0^t\langle\A \theta,
\theta_t\rangle=-\int_0^t\int_\Omega \theta
h'(\theta)\theta_t\chi_t+\int_0^t\int_\Omega \chi_t^2\theta_t.
\label{t.1}
\end{equation}
Owing to \eqref{h2}-\eqref{h3} we obtain
\begin{align}
&c_s\norm{\theta_t}^2_{L^2(0,t,H)}\leq \int_0^t\int_\Omega
\theta_t^2
(1+\chi\theta h''(\theta))\label{t.2}\\
\intertext{and} &\abs{\int_0^t\int_\Omega \theta h'(\theta) \theta_t
\chi_t}+\abs{\int_0^t\int_\Omega \chi_t^2\theta_t}\leq
c_h\int_0^t\int_\Omega
\abs{\chi_t}\abs{\theta_t}+\int_0^t\int_\Omega
\abs{\chi_t}^2\abs{\theta_t}\notag\\
&\leq C_{16}\int_0^t\left(\norm{\chi_t}_H+\norm{\chi_t}_V^2\right)\norm{\theta_t}_H\notag\\
&\leq\frac{c_s}{2}\norm{\theta_t}^2_{L^2(0,t,H)}+C_{17}\int_0^t
\norm{\chi_t}^2_H+C_{18}\int_0^t\norm{\chi_t}^2_V
\norm{\chi_t}^2_V,\label{t.3}
\end{align}
where we have exploited the Young and H\"{o}lder inequalities as
well as the continuous embedding $V\subset L^4(\Omega)$. By
integrating by parts in time, thanks to \eqref{teta0} it follows
that
\begin{equation} \label{t.4}
\frac{c_s}{2} \norm{\theta_t}^2_{L^2(0,t,H)}+\frac{1}{2}\norm{\nabla
\theta (t)}^2_H \leq C_{19}+C_{17}\int_0^t
\norm{\chi_t}^2_H+C_{18}\int_0^t\norm{\chi_t}^2_V \norm{\chi_t}^2_V,
\end{equation}
where for the moment we just know that $\norm{\chi_t}^2_V\in
L^1(0,T)$ (see \eqref{chiregularity}). However, we are going to
combine estimate \eqref{t.4} with another estimate which will give
us more information on $\norm{\chi_t}^2_V$. Indeed, we consider
\eqref{2} and differentiate it with respect to time thus obtaining
\begin{equation} \label{c.1}
\chi_{tt}+\nu\A\chi_{tt}+\A\chi_t+\xi_t=h'(\theta)\theta_t-\frac{1}{p}p_t.
\end{equation}
Then, we test \eqref{c.1} by $\chi_t$. After some integrations by
parts in time and owing to \eqref{h1}, we can write
\begin{align} \label{c.2}
&\frac{1}{2}\int_\Omega\left(\abs{\chi_t}^2+\nu\abs{\nabla\chi_t}^2\right)(t)+\int_0^t\int_\Omega
\abs{\nabla\chi_t}^2+\int_0^t\langle\xi_t,\chi_t\rangle\notag\\
&\leq C_{20}
\norm{\chi_t(0)}^2_V+c_h\int_0^t\int_\Omega\abs{\theta_t}\abs{\chi_t}+\int_0^t\int_\Omega\abs{\frac{1}{p}}\abs{p_t}\abs{\chi_t}.
\end{align}
Let us note that since $\chi_0\in D(\partial_{V,V'}J)$, there exists
$\xi_0\in\partial_{V,V'}J(\chi_0)$. Thus, we can introduce $\chi'_0$
as the initial value of the time derivative of $\chi$ by defining it
as the solution of the following elliptic equation
\begin{equation}
\chi'_0+\nu\A\chi'_0=-\A\chi_0-\xi_0+h(\theta_0)-\log p_0,
\end{equation}
where the right hand side is known in $V'$ thanks to
\eqref{teta0}-\eqref{p0}, \eqref{h1} and \eqref{p^-10}. Then, it
easily follows that $\chi'_0\in V$. Next, we observe that
monotonicity arguments yield
\begin{equation} \label{monotonicity}
\int_0^t\langle\xi_t,\chi_t\rangle\geq 0.
\end{equation}
By applying H\"{o}lder's and Young's inequalities, we can estimate
the first integral on the right hand side of \eqref{c.2} as follows
\begin{equation} \label{c.3}
c_h\int_0^t\int_\Omega \abs{\theta_t}\abs{\chi_t}\leq
\frac{c_s}{4}\norm{\theta_t}^2_{L^2(0,t,H)}+C_{21}\int_0^t
\norm{\chi_t}^2_H.
\end{equation}
Analogously, the last integral on the right hand side of \eqref{c.2}
can be handled as follows
\begin{equation} \label{c.4}
\int_0^t \int_\Omega \abs{\frac{1}{p}}\abs{p_t}\abs{\chi_t}\leq
\int_0^t
\norm{\frac{1}{p}}_{L^4(\Omega)}\norm{p_t}_H\norm{\chi_t}_V,
\end{equation}
where, owing to Theorems \ref{teorexistence} and
\ref{teorregularitypressure},
$(\norm{1/p}_{L^4(\Omega)}\norm{p_t}_H)\in L^1(0,T)$. Thanks to
\eqref{monotonicity}-\eqref{c.4}, \eqref{c.2} yields
\begin{align} \label{c.5}
&\frac{1}{2}\int_\Omega\left(\abs{\chi_t}^2+\nu\abs{\nabla\chi_t}^2\right)(t)+\int_0^t\int_\Omega
\abs{\nabla\chi_t}^2\notag\\
&\leq
C_{22}+\frac{c_s}{4}\norm{\theta_t}^2_{L^2(0,t,H)}+C_{21}\int_0^t
\norm{\chi_t}^2_H+\int_0^t
\norm{\frac{1}{p}}_{L^4(\Omega)}\norm{p_t}_H\norm{\chi_t}_V.
\end{align}
Finally, by combining \eqref{t.4} and \eqref{c.5}, we obtain
\begin{align} \label{combination}
&\frac{c_s}{4}\norm{\theta_t}^2_{L^2(0,t,H)}+\frac{1}{2}\norm{\nabla
\theta(t)}^2_H+\frac{1}{2}\int_\Omega\left(\abs{\chi_t}^2+\nu\abs{\nabla\chi_t}^2\right)(t)+\norm{\nabla\chi_t}^2_{L^2(0,t,H)}\notag\\
&\leq C_{23}+C_{18}\int_0^t\norm{\chi_t}^2_V
\norm{\chi_t}^2_V+C_{24}\int_0^t \norm{\chi_t}^2_H+\int_0^t
\norm{\frac{1}{p}}_{L^4(\Omega)}\norm{p_t}_H\norm{\chi_t}_V.
\end{align}
Now, we can apply the Gronwall lemma \cite[Theorem 2.1]{ba} and
owing to Theorem \ref{teorexistence} (see \eqref{tetaregularity} and
\eqref{chiregularity}) we get
\begin{equation} \label{tc}
\norm{\theta}_{H^1(0,T,H)\cap L^\infty
(0,T,V)}+\norm{\chi}_{W^{1,\infty}(0,T,V)}\leq c,
\end{equation}
so that Theorem \ref{teorfurtherregularities} is completely proved.

Finally, we aim to prove Theorem \ref{teorpositivitytemperature},
i.e. to establish the positivity of the temperature. To this
purpose, we deal with \eqref{t.0} and formally test it by
$-\theta^{-1}$. After integrating over $(0,t)$, for the first term
we get
\begin{align}
-\int_0^t\int_\Omega \theta_t\theta^{-1}&=-\int_0^t\int_\Omega
\frac{d}{dt}(\log \theta)\notag\\
&=\int_\Omega (\log\theta)^-(t)-\int_\Omega(\log
\theta)^+(t)+\int_\Omega \log \theta_0,\label{pos.1}
\end{align}
where $(\log\theta)^-$ and $(\log\theta)^+$ denote the negative and
positive parts of the function $\log\theta$, respectively. Hence, we
have
\begin{align}
&\int_\Omega \abs{\log \theta}(t)-\int_0^t\int_\Omega
h''(\theta)\theta_t \chi+\int_0^t\int_\Omega\abs{\nabla\log
\theta}^2+\int_0^t\int_\Omega
\frac{\chi^2_t}{\theta}\notag\\
&=2\int_\Omega(\log \theta)^+(t)-\int_\Omega \log
\theta_0+\int_0^t\int_\Omega h'(\theta)\chi_t.\label{pos.2}
\end{align}
Now, we deal with the second integral on the left hand side of
\eqref{pos.2} and, after integrating by parts in time, we can write
\begin{align}
-\int_0^t\int_\Omega
h''(\theta)\theta_t\chi&=-\int_0^t\int_\Omega\frac{\partial}{\partial t}(h'(\theta))\chi\notag\\
&=-\int_\Omega h'(\theta(t))\chi(t)+\int_\Omega
h'(\theta_0)\chi_0+\int_0^t\int_\Omega
h'(\theta)\chi_t,\label{pos.3}
\end{align}
and then
\begin{align}
&\int_\Omega \abs{\log \theta}(t)+\int_0^t\int_\Omega\abs{\nabla\log
\theta}^2+\int_0^t\int_\Omega
\frac{\chi^2_t}{\theta}\notag\\
&=\int_\Omega h'(\theta(t))\chi(t)+2\int_\Omega(\log
\theta)^+(t)-\int_\Omega \log \theta_0-\int_\Omega
h'(\theta_0)\chi_0.\label{pos.4}
\end{align}
Let us observe that well-known properties of the logarithm function
and \eqref{tetaregularity} yield
\begin{align}
\norm{(\log\theta)^+}_{L^\infty(0,T,L^1(\Omega))}&\leq \sup _{0\leq t\leq T}\int_{\Omega\cap\{\theta\geq1\}}\abs{\log\theta(t)}\notag\\
&\leq \sup _{0\leq t\leq
T}\int_{\Omega\cap\{\theta\geq1\}}\abs{\theta}^2\leq
\norm{\theta}^2_{L^\infty(0,T,H)}\leq c.\label{pos.5}
\end{align}
Besides, we note that the first integral in the right hand side of
\eqref{pos.4} is bounded thanks to \eqref{h1} and
\eqref{chiregularity}. By virtue of \eqref{chi0}, \eqref{logteta0}
and \eqref{pos.5} we can eventually infer that
\begin{equation} \label{pos.6}
\norm{\log \theta(t)}_{L^1(\Omega)}+\norm{\nabla\log
\theta}^2_{L^2(0,T,H)}+\norm{\frac{\chi_t}{\sqrt{\theta}}}^2_{L^2(0,T,H)}\leq
c.
\end{equation}
As a consequence of \eqref{pos.5}, \eqref{pos.6} and exploiting the
Poincar\'e-Wirtinger ine\-qua\-li\-ty, we can deduce
\begin{equation} \label{pos.final}
\norm{\log\theta}_{L^\infty(0,T,L^1(\Omega))\cap
L^2(0,T,V)}+\norm{\frac{\chi_t}{\sqrt{\theta}}}^2_{L^2(0,T,H)}\leq
c.
\end{equation}

\section{Appendix}
In this appendix, we present two auxiliary results among those
exploited in our proofs. As first, let us detail and prove a
convergence result we have applied in order to pass to the limit in
the time discretization scheme.
\begin{propos} \label{passaggioallimite}
Let $\mathcal{T_\tau}$ be the translation operator defined in
\eqref{translation}. If the sequence $\{v_\tau\}$ satisfies
\begin{equation} \label{ipotesipassaggioallimite}
v_\tau\rightarrow v \quad\text{in}\quad L^2(-T,T,H),
\end{equation}
as $\tau\searrow0$, then the following convergence
\begin{equation} \label{tesipassaggioallimite}
\int_0^T\norm{\mathcal{T_\tau}v_\tau(t)- v(t)}^2_H dt
\xrightarrow[\tau\searrow 0]{} 0
\end{equation}
holds.
\end{propos}
\begin{pr}
Here, we give just an outline of the proof. A first step consists in
proving that the following convergence
\begin{equation} \label{step1}
\int_0^T\norm{\mathcal{T_\tau}w(t)- w(t)}^2_H dt
\xrightarrow[\tau\searrow 0]{} 0
\end{equation}
holds for any $w\in C^0([-T,T],H)$\footnote{By $C^0([-T,T],H)$ we
denote the space of continuous functions from $[-T,T]$ into $H$
equipped with the $L^\infty(-T,T,H)$ norm.}. To this aim, let us
note that, if $w\in C^0([-T,T],H)$, then \eqref{step1} easily
follows by suitably applying the Lebesgue's dominated convergence
theorem. Now, since $C^0([-T,T], H)$ is dense in $L^2(-T,T,H)$, if
$v\in L^2(-T,T,H)$, there exists an approximating sequence
${\{v_n\}}_{n \in \mathbb{N}}\subset C^0([-T,T], H)$ such that
\begin{equation} \label{step2}
v_n\rightarrow v\quad\text{in}\quad L^2(-T,T,H)
\end{equation}
as $n\rightarrow+\infty$. In particular, for any $\varepsilon>0$,
there exists $n_\varepsilon\in \mathbb{N}$ such that the inequality
\begin{equation}\label{limite1}
\norm{v_n-v}^2_{L^2(-T,T,H)}\leq \varepsilon
\end{equation}
holds for all $ n\geq n_\varepsilon$. Moreover, thanks to
\eqref{translation}, from \eqref{limite1} we can easily deduce that
\begin{equation} \label{step3}
\int_0^T\norm{\mathcal{T_\tau}v_n(t)-\mathcal{T_\tau}v(t)}^2_H
dt\leq\norm{v_n-v}^2_{L^2(-T,T,H)}\leq\varepsilon
\end{equation}
for all $n\geq n_\varepsilon$ and for all $0<\tau<T$. Now, we fix
$n=n_\varepsilon$ and point out that since $v_n\in C^0([-T,T], H)$
then \eqref{step1} holds for $w=v_n$. In particular, there exists
$\tau_\varepsilon$ such that
\begin{equation}
\int_0^T\norm{\mathcal{T_\tau}v_n(t)- v_n(t)}^2_H dt\leq\varepsilon
\end{equation}
for all $\tau\leq\tau_\varepsilon$. This allows us to split the
difference $(\mathcal{T_\tau}v-v)$ into three terms, all tending to
zero in $L^2(0,T,H)$. Namely, we can write
\begin{align} \label{split1}
\mathcal{T_\tau}v-v&=\mathcal{T_\tau}v-\mathcal{T_\tau}
v_n\notag\\
&+\mathcal{T_\tau}v_n-v_n\notag\\
&+v_n-v.
\end{align}
By means of the above arguments, we can infer that the following
inequality
\begin{equation}
\int_0^T\norm{\mathcal{T_\tau}v(t)-v(t)}^2_H dt\leq3\varepsilon
\end{equation}
holds for all $\tau\leq \tau_\varepsilon$, whence
\begin{equation}\label{ipotesinew}
\int_0^T\norm{\mathcal{T_\tau}v(t)-v(t)}^2_H dt
\xrightarrow[\tau\searrow 0]{} 0.
\end{equation}
Finally, we consider the difference between $\mathcal{T_\tau}v_\tau$
and $v$. It is convenient to express it as follows
\begin{align} \label{split2}
\mathcal{T_\tau}v_\tau-v&=\mathcal{T_\tau}v_\tau-\mathcal{T_\tau}
v\notag\\
&+\mathcal{T_\tau}v-v.
\end{align}
Thanks to \eqref{ipotesipassaggioallimite} and to the definition of
$\mathcal{T_\tau}$, for the first addend on the right hand side we
can write
\begin{equation}
\int_0^T\norm{\mathcal{T_\tau}v_\tau(t)-\mathcal{T_\tau}
v(t)}^2_H\leq
\norm{v_\tau-v}^2_{L^2(-T,T,H)}\xrightarrow[\tau\searrow 0]{} 0,
\end{equation}
for all $0<\tau<T$. While the second addend has already been
discussed in \eqref{ipotesinew}. Hence, we can infer that
\begin{equation}
\int_0^T\norm{\mathcal{T_\tau}v_\tau(t)- v(t)}^2_H dt
\xrightarrow[\tau\searrow 0]{} 0
\end{equation}
and the proof is complete.
\end{pr}
Finally, we present a fairly standard result concerning the sum of
maximal monotone operators. Let us recall some notations previously
introduced. We denote by $H$ the space $L^2(\Omega)$, with $\Omega$
a bounded domain included in $\R^3$, and by $\A$ the abstract
operator prescribed in \eqref{A}. Next, by $j$ we denote a proper,
lower semicontinuous and convex function on $H$, and by $\partial_H
j$ the subdifferential of $j$. Then the following result holds.
\begin{teor} \label{teorsuriectivitymaximalmonotoneoperators}
Let $f\in H$. Then, there exists a unique $\chi\in H^2(\Omega)\cap
D(j)$ such that
\begin{equation} \label{equationmaxmonop}
\chi+\A \chi+\partial_H j (\chi) \ni f.
\end{equation}
\end{teor}
The above result can be easily proved by approximating the
subdifferential operator in \eqref{equationmaxmonop} by its Yosida
regularization $(\partial_H j) _\lambda$. We remind that
$(\partial_H j) _\lambda$ is maximal monotone and Lipschitz
continuous. Hence, we may invoke well known results on the sum of
maximal monotone operators (cf. e.g.  \cite{bre}) and get the
existence of a solution $\chi_\lambda\in H^2(\Omega)$ solving the
approximated equation. Next, we could perform standard a priori
estimates, independent of $\lambda$, on the approximated system and
get the limit as $\lambda\searrow0$ by compactness. Finally,
uniqueness could be deduced via standard contradiction arguments
owing to the monotonicity of the subdifferential operator.

\section*{Acknowledgements}
I am deeply grateful to Professor Pierluigi Colli for his
encouragement, support and cooperation to this work. My teachers of
Mathematics at the University of Pavia and the kind hospitality and
stimulating atmosphere of the Mathematics Department in Pavia are
gratefully acknowledged.

\end{document}